\documentclass[onefignum,onetabnum]{siamonline250106}

\usepackage[ansinew]{inputenc}
\usepackage{amsmath, amsfonts} 
\usepackage{amssymb}
\usepackage[english]{babel}
\usepackage[T1]{fontenc}
\usepackage{lmodern}
\usepackage{graphicx}
\usepackage{subfigure}
\usepackage{lipsum}
\usepackage{graphicx}
\usepackage{verbatim}
\usepackage{array}
\usepackage{float}
\usepackage{subfigure}
\usepackage{color}
\usepackage{algorithm, algorithmicx, algpseudocode}

\newcommand{\Hc}{\mathcal{H}}
\ifpdf
  \DeclareGraphicsExtensions{.eps,.pdf,.png,.jpg}
\else
  \DeclareGraphicsExtensions{.eps}
\fi

\usepackage{enumitem}
\setlist[enumerate]{leftmargin=.5in}
\setlist[itemize]{leftmargin=.5in}

\newsiamremark{remark}{Remark}
\newsiamremark{example}{Example}
\newsiamremark{hypothesis}{Hypothesis}
\crefname{hypothesis}{Hypothesis}{Hypotheses}
\newsiamthm{claim}{Claim}

\headers{Analytic EDMD}{A. Mauroy and I.Mezi\'c}

\title{Analytic Extended Dynamic Mode Decomposition\thanks{Submitted to the editors DATE.
}}

\author{A. Mauroy\thanks{Department of Mathematics and Namur Institute for Complex Systems, University of Namur, Belgium 
  (\email{alexandre.mauroy@unamur.be}).}
\and I. Mezi\'c\thanks{Department of Mechanical Engineering, University of California Santa Barbara, US
  (\email{mezic@ucsb.edu}).}
}

\usepackage{amsopn}

\ifpdf
\hypersetup{
  pdftitle={Analytic EDMD},
  pdfauthor={A. Mauroy and I. Mezi\'c}
}
\fi

\begin{document}

\maketitle

\begin{abstract}
We develop a novel EDMD-type algorithm that captures the spectrum of the Koopman operator defined on a reproducing kernel Hilbert space of analytic functions. This method, which we call analytic EDMD, relies on an orthogonal projection on polynomial subspaces, which is equivalent to a data-driven Taylor approximation. In the case of dynamics with a hyperbolic equilibrium, analytic EDMD demonstrates excellent performance to capture the lattice-structured Koopman spectrum based on the eigenvalues of the linearized system at the equilibrium. Moreover, it yields the Taylor approximation of associated principal eigenfunctions. Since the method preserves the triangular structure of the operator, it does not suffer from spectral pollution and, moreover, arbitrary accuracy on the spectrum can be reached with a fixed finite dimension of the approximation and with a (possibly non-uniform) sampling over an arbitrary set of nonzero measure. The performance of analytic EDMD is illustrated with numerical examples and is assessed through a comparative study with related methods. Finally, the method is complemented with theoretical results, proving strong convergence of the eigenfunctions and providing error bounds on the spectrum estimation.
\end{abstract}

\begin{keywords}
Dynamic Mode Decomposition, Koopman Operator, Kernel Methods, Taylor approximation, Reproducing kernel Hilbert space
\end{keywords}

\begin{MSCcodes}
34L16, 37M10, 41A58, 47A58, 46E20, 46E22
\end{MSCcodes}

\section{Introduction}

For several decades, the Koopman (or composition) operator has been thoroughly studied from a theoretical point of view on spaces of analytic functions \cite{Cowen1995composition}. In the case of equilibrium dynamics generated by analytic discrete-time maps or vector fields, recent works have investigated the properties of the Koopman operator defined on different functional spaces: the Hardy space on the polydisk \cite{Zagabe_SICON} (see also \cite{Rosenfeld_Hardy} for general properties) and the (modulated) Fock space of entire analytic functions \cite{Mezic2020}. In this case, the spectrum of the Koopman operator is known to possess a lattice structure including the eigenvalues of the linearized system at the equilibrium (see e.g. \cite{Maccluer1984spectra,Mezic2020}) and the corresponding principal eigenfunctions reveal meaningful properties of the dynamics (e.g. global stability \cite{Mauroy_Mezic_stability}, isostables \cite{MMM_isostables}). Interestingly, it is well-known that a finite-dimensional Koopman representation obtained via a Taylor-based truncation of analytic functions yields the \emph{exact} Koopman eigenvalues (see e.g. \cite[Chapter 1]{book_Koopman}) and the \emph{exact} Taylor approximation of the associated eigenfunctions, a property which is explained by the block-triangular structure of the Koopman operator approximation \cite{Zagabe_SICON} (see also the Schur decomposition discussed in \cite{Drmac_Schur}).

In contrast, the standard extended dynamic mode decomposition (EDMD) method \cite{Rowley_EDMD} relying on an $L^2$ (Galerkin) projection does not preserve the block-triangular structure and is in particular affected by spectral pollution \cite{Colbrook2024, Drmac2018, Korda_convergence}. Also, in \cite[Theorem 3]{Mezic2022}, it is shown that the coefficients in the $L^2$ projection-based approximation of the Koopman operator are equal to ergodic averages over the attractor of the dynamical system. This implies that, in the single-trajectory case, the approximation matrix obtained in the limit of large data (infinite time series) can ``see'' only the attractor dynamics and thus cannot capture the transient dynamics (off-attractor). In other words, a single trajectory-based $L^2$ approximation can only reveal on-attractor eigenvalues. Off-attractor eigenvalues are better captured with a set of independent trajectories (see e.g. \cite{Lucarini} and other set-oriented methods \cite{Dellnitz2001}). But yet, in this setting, spectral convergence results for EDMD applied to analytic maps have only been proven recently under very specific conditions (e.g. one-dimensional systems) \cite{Slipantschuk_interval}.

For the above-mentioned reasons, it appears that EDMD is not appropriate to reveal spectral properties of analytic transient dynamics. Instead, a data-driven Taylor projection is much more desirable, and actually happens to be the center of recent attention. In \cite{Ishikawa_JetEDMD}, this idea has been approached by combining the standard EDMD method with a finite-section truncation, a method referred to as JetEDMD which was complemented with convergence results and shown to yield improved performance compared to EDMD. This method also provides generalized eigenfunctions related to an extended Koopman operator defined in a rigged Hilbert space, although their theoretical interpretation is still not clear, as acknowledged by the authors. Moreover, the work \cite{EDMD_analytic_circle} computes orthogonal projections in the Hardy-Hilbert space on a complex annulus to design a novel EDMD-type algorithm that captures the Koopman spectrum of analytic expanding circle maps. This line of research exploits the fact that weighted monomials form an orthonormal basis in weighted Hardy spaces of analytic functions, so that the Taylor projection is an orthogonal projection. However, inner products are typically expressed as integrals over the complex circle (e.g. in the Hardy space) or over the complex plane (e.g. in the Fock space). And computing such integrals is not convenient from a practical (data-driven) point of view, as dynamics of interest typically evolve on a real state space and therefore cannot generate complex-valued data. In the present work, we circumvent this issue by exploiting the properties of reproducing kernel Hilbert spaces (RKHS). Using a kernel associated with a space of analytic functions, we derive a simple data-driven Taylor projection based on real-valued data sets and we leverage this technique to design the so-called \emph{analytic EDMD} algorithm. This method is reminiscent of kernel ridge regression \cite{Klus_kernel, Kostic} and loosely connected to occupation kernel-based methods developed in weighted Hardy spaces \cite{Rosenfeld_singular_DMD, Rosenfeld_DMD_Liouville}. Analytic EDMD is shown to demonstrate excellent performance to capture spectral properties of analytic dynamics with a hyperbolic equilibrium, in particular revealing the lattice-structured Koopman spectrum (including the eigenvalues of the linearized dynamics) and the principal eigenfunctions, even with data taken far from the equilibrium. Since the method preserves the block-triangular structure of the operator, it does not suffer from spectral pollution. Moreover, theoretical results show strong convergence of Koopman eigenfunctions and yield error bounds for eigenvalue estimation which only depend on sample complexity. It follows that arbitrary accuracy on the spectrum can be obtained with a fixed finite dimension of the approximation and, thanks to the properties of analytic functions, with a (possibly non-uniform) sampling over a set of nonzero Lebesgue measure.

The rest of the paper is organized as follows. In \cref{sec:Taylo_proj}, we derive the data-driven Taylor projection technique that is leveraged in \cref{sec:analytic_EDMD} to develop the analytic EDMD method. Illustrative examples are provided in \cref{sec:num_ex} and numerical analysis including a comparative study is performed in \cref{sec:num_analysis}. Theoretical convergence results and error bounds are obtained in \cref{sec:errors}. Finally, concluding remarks and perspectives are given in \cref{sec:conclu}.

\section{Data-driven Taylor projection}
\label{sec:Taylo_proj}

In this section, we leverage the RKHS property of analytic spaces to obtain a data-driven approximation of inner products. This yields a data-driven approximation of the orthogonal Taylor projection, which we will use in our analytic EDMD method.

Let us consider an infinite-dimensional RKHS $\Hc$ over the compact set $X \subset \mathbb{R}^n$ associated with the continuous kernel $k:X\times X \to \mathbb{R}$. We assume that the (weighted) monomials $\{e_i(\mathbf{x})\}_{i=1}^\infty=\{\beta_\alpha \mathbf{x}^\alpha\}_{\alpha \in \mathbb{N}^n}$ form an orthonormal basis in $\Hc$, i.e. $\langle e_i, e_j \rangle_\Hc=\delta_{ij}$. This is typically the case for RKHSs associated with Taylor-type kernels of the form $k(\mathbf{x},\mathbf{y})=f(\mathbf{x}^T\mathbf{y})$, where $f:X \to \mathbb{R}$ is an analytic function with non-negative Taylor coefficients (see e.g. \cite{SVD_book}).

Let $m$ be a Borel probability measure on $X$ and $T:L^2(X,m) \to L^2(X,m)$ be the kernel integral operator 
\begin{equation*}
(Tf)(\mathbf{x}) = \int_X k(\mathbf{x},\mathbf{y}) \, f(\mathbf{y}) \, dm(\mathbf{y}).
\end{equation*}
It follows from Mercer's theorem that $T^{1/2}$ defines an isometry\footnote{We suppose that the kernel integral operator is injective. This is equivalent to $\Hc$ being dense in $L^2(X,m)$ (\cite[Theorem 4.26]{SVD_book}), which holds in the case of an analytic RKHS.} between $\Hc$ and $L^2(X,m)$ (see e.g. \cite[Chapter 11]{RKHS_Paulsen} and \cite[Chapter 4]{SVD_book}), so that the inner product of $\Hc$ can be rewritten as
\begin{equation}
\label{eq:inner_prod}
\langle f,g \rangle_\Hc = \langle T^{-1/2} f,T^{-1/2} g \rangle_{L^2(X,m)}.
\end{equation}
For a given set of sample points $\{\mathbf{x}_k\}_{k=1}^M \subseteq X$ randomly distributed according to the measure $m$, we have the empirical estimate
\begin{equation*}
(Tf)(\mathbf{x}_i) \approx \frac{1}{M} \sum_{j=1}^M k(\mathbf{x}_i,\mathbf{x}_j) \, f(x_j)
\end{equation*}
or equivalently
\begin{equation*}
\begin{pmatrix}
(Tf)(\mathbf{x}_1) \\ \vdots \\ (Tf)(\mathbf{x}_M)
\end{pmatrix}
\approx \frac{1}{M} \begin{pmatrix}
k(\mathbf{x}_1,\mathbf{x}_1) & \cdots & k(\mathbf{x}_1,\mathbf{x}_M) \\
\vdots & \ddots & \vdots \\
k(\mathbf{x}_1,\mathbf{x}_M) & \cdots & k(\mathbf{x}_M,\mathbf{x}_M)
\end{pmatrix} 
\begin{pmatrix}
f(\mathbf{x}_1) \\ \vdots \\ f(\mathbf{x}_M)
\end{pmatrix},
\end{equation*}
where
\begin{equation}
\label{eq:kernel_mat}
 \begin{pmatrix}
k(\mathbf{x}_1,\mathbf{x}_1) & \cdots & k(\mathbf{x}_1,\mathbf{x}_M) \\
\vdots & \ddots & \vdots \\
k(\mathbf{x}_1,\mathbf{x}_M) & \cdots & k(\mathbf{x}_M,\mathbf{x}_M)
\end{pmatrix} \triangleq \mathbf{G} 
\end{equation}
is the symmetric (Gram) kernel matrix. It follows that the inner product \eqref{eq:inner_prod} can be approximated by
\begin{equation*}
\langle f,g \rangle_\Hc = \langle T^{-1/2} f,T^{-1/2} g \rangle_{L^2(X,m)} \approx \frac{1}{M}  \left(\left(\frac{1}{M}\mathbf{G}\right)^{-1/2} \mathbf{f}\right)^T \left(\left(\frac{1}{M} \mathbf{G}\right)^{-1/2} \mathbf{g}\right)  = \mathbf{f}^T  \mathbf{G}^{-1} \mathbf{g}
\end{equation*}
with $\mathbf{f} = (f(\mathbf{x}_1)\, \cdots \, f(\mathbf{x}_M))^T$ and $\mathbf{g} = (g(\mathbf{x}_1)\, \cdots \, g(\mathbf{x}_M))^T$.

\newpage
\begin{remark}[Interpretation of the inner product approximation\footnote{This remark is from a private communication with Isao Ishikawa.}]
\label{rem:inner_product}
    The proposed data-driven approximation of $\langle f,g \rangle_\Hc$ can be interpreted as the inner product $\langle \Pi f,g \rangle_\Hc=\langle f,\Pi g \rangle_\Hc$, where $\Pi$ is the orthogonal projection onto the span of reproducing functions $k_{\mathbf{x}_k}=k(\mathbf{x}_k, \cdot)$. Indeed, it is well-known that (see e.g. \cite[Chapter 3]{RKHS_Paulsen})
 \begin{equation*}
    \Pi f = \mathbf{f}^T \mathbf{G}^{-1} \begin{pmatrix}
        k_{\mathbf{x}_1} \\ \vdots \\ k_{\mathbf{x}_M}
    \end{pmatrix}
\end{equation*}   
and, using the reproducing property $\langle k_{\mathbf{x}_k}, g \rangle_\Hc = g(\mathbf{x}_k)$, we obtain
\begin{equation}
\label{eq:prop_inner_prod}
    \langle \Pi f,g \rangle_\Hc = \mathbf{f}^T \mathbf{G}^{-1} \begin{pmatrix}
        \langle k_{\mathbf{x}_1},g \rangle_\Hc \\ \vdots \\ \langle k_{\mathbf{x}_M}, g \rangle_\Hc \end{pmatrix}
        = \mathbf{f}^T \mathbf{G}^{-1} \mathbf{g} .
\end{equation}
\hfill $\diamond$
\end{remark}

Next, suppose that $\Hc_N \subset \Hc$ is a subspace spanned by an orthonormal basis of $N$ weighted monomials $\{e_i\}_{i=1}^N$. Then, the (Taylor) orthogonal projection $P:\Hc \to \Hc_N$ can be approximated from the data by
\begin{equation}
\label{eq:Taylor_proj}
 Pf = \sum_{i=1}^N \langle f,e_i \rangle_\Hc \, e_i \approx \sum_{i=1}^N \left(\mathbf{f}^T  \mathbf{G}^{-1} \mathbf{e}_i\right) \, e_i
\end{equation}
with $\mathbf{f} = (f(\mathbf{x}_1)\, \cdots \, f(\mathbf{x}_M))^T$ and $\mathbf{e}_i = (e_i(\mathbf{x}_1)\, \cdots \, e_i(\mathbf{x}_M))^T$.

\begin{example}
    We approximate the Taylor coefficients (around $x=0$) of the function $f(x)=1+\log(x)$, using $10$ data points $(x_k,f(x_k))$ drawn from a uniform distribution over $[-1,1]$. The results shown in \cref{tab1} are obtained with the Szeg\"o kernel $k(x,y)= 1/(1-x y)$ associated with the Hardy space on the polydisc, and with a basis of monomials up to degree $5$. They are compared to the results obtained with a polynomial regression based on a discrete $L^2$ projection.
\begin{table}[h]
    \centering
    \begin{tabular}{ccccccc}
        degree & 0 & 1 & 2 & 3 & 4 & 5 \\
        \hline
        exact & 0 & 1 & -0.5 & 0.333 & -0.25 & 0.2 \\
        Taylor projection & -1.04e-6 & 1.000 & -0.500 & 0.334 & -0.245 & 0.202 \\
        $L^2$ projection & -0.01 & 1.010 & -0.416 & 0.262 & -0.555 & 0.423 \\
    \end{tabular}
    \caption{Data-driven computation of the Taylor coefficients of $f(x)=1+\log(x)$.}
    \label{tab1}
\end{table}
    
\end{example}

\section{Analytic EDMD}
\label{sec:analytic_EDMD}

We are now in a position to develop the analytic EDMD method, based on a data-driven finite section approximation of the Koopman operator \mbox{$K:\Hc \to \Hc$}, $Kf= f \circ \varphi$, where $\varphi$ is an analytic map $\varphi:X \to X$ with a unique equilibrium in $X$. Note that $\varphi$ can be interpreted as a discrete-time map or as the flow of a continuous-time system evaluated at a fixed sampling time $\Delta t$. In the latter case, an analytic flow is typically generated by an analytic vector field that possesses a hyperbolic equilibrium with non-resonant Jacobian eigenvalues (Poincar\'e linearization theorem, see e.g. \cite{Gaspard}). We also make the standing assumption that the RKHS $\Hc$ is invariant under the action of the Koopman operator, and in particular $K e_i \in \Hc$ for all monomials $e_i$. The RKHS should therefore be infinite-dimensional, even for a generic polynomial map $\varphi$ or polynomial vector field.

If the dynamics admit an equilibrium at the origin (without loss of generality), it is well-known that the Koopman operator satisfies $\langle e_i, K e_j \rangle_\Hc = 0$ if the total degree of the monomial $e_i$ is smaller than the total degree of $e_j$, so that its matrix representation in a monomial basis is lower block-triangular (provided that the monomials are sorted by increasing total degree). This is explained by the fact that the composition of a monomial with an analytic map yields a Taylor expansion of higher degree. Note that a similar property also holds with the Koopman infinitesimal generator (see e.g. \cite{Zagabe_SICON}). Equivalently, the Koopman operator satisfies $P K P = P K$ if $P:\Hc \to \Hc_N$ is the orthogonal Taylor projection onto a subspace of polynomials up to a given total total degree. If $\mu \in \sigma(K)$ is an eigenvalue of $K$ associated with the eigenfunction $\phi_\mu$, it follows that $P K P \phi_\mu = P K \phi_\mu = \mu P \phi_\mu$, so that $\sigma(P K P ) \subset \sigma(K)$ and $P \phi_\mu$ is an eigenfunction of $PKP$. This shows that the spectral properties of the Koopman operator are recovered from its approximation $PKP$, an observation which motivates the use of the finite section method obtained via the data-driven Taylor projection.

\subsection{Finite-dimensional approximation}
\label{sec:finite_dim_approx}

The entries of the Koopman matrix approximation $\mathbf{\hat{K}}$ obtained with the finite section method are such that $P Ke_j = \sum_{i=1}^N \mathbf{K}_{ij} \, e_i$. Using \cref{eq:Taylor_proj}, we have
\begin{equation*}
P K e_j = \sum_{i=1}^N \langle K e_j,e_i \rangle_\Hc \, e_i \approx \sum_{i=1}^N \left((\mathbf{e}_j')^T  \mathbf{G}^{-1} \mathbf{e}_i \right) \, e_i
\end{equation*}
with $(\mathbf{e}_j')^T=((Ke_j)(\mathbf{x}_1)\, \cdots \, (Ke_j)(\mathbf{x}_M))^T$. Then the entries $\mathbf{K}_{ij}=\langle K e_j,e_i \rangle_\Hc$ are approximated by $\hat{\mathbf{K}}_{ij} = (\mathbf{e}_i)^T  \mathbf{G}^{-1} \mathbf{e}_j'$. Assuming that we are given $M$ data pairs $(\mathbf{x}_k,\mathbf{y}_k)_{k=1}^M=(\mathbf{x}_k,\varphi(\mathbf{x}_k))_{k=1}^M$, we can rewrite \mbox{$(\mathbf{e}_j')^T=(e_j(\mathbf{y}_1)\, \cdots \, e_j(\mathbf{y}_M))^T$} and we obtain
\begin{equation}
\label{eq:analytic_EDMD}
\hat{\mathbf{K}} = \mathbf{X}^T \mathbf{G}^{-1} \mathbf{Y} 
\end{equation}
with the data matrices
\begin{equation}
\label{eq:data_mat}
\mathbf{X} = \begin{pmatrix}
e_1(\mathbf{x}_1) & \cdots & e_N(\mathbf{x}_1) \\
\vdots & \ddots & \vdots \\
e_1(\mathbf{x}_M) & \cdots & e_N(\mathbf{x}_M)
\end{pmatrix} 
\qquad
\mathbf{Y} = \begin{pmatrix}
e_1(\mathbf{y}_1) & \cdots & e_N(\mathbf{y}_1) \\
\vdots & \ddots & \vdots \\
e_1(\mathbf{y}_M) & \cdots & e_N(\mathbf{y}_M)
\end{pmatrix}. 
\end{equation}
Let us assume that the monomials $e_i(\mathbf{x})=\beta_{\overline{\alpha}(i)} \mathbf{x}^{\overline{\alpha}(j)}$ are sorted according to the (possibly lexicographic) order defined by the map $\overline{\alpha}:\mathbb{N} \to \mathbb{N}^n$ such that $i_1<i_2$ if $|\overline{\alpha}(i_1)|<|\overline{\alpha}(i_2)|$, where $|\overline{\alpha}(i)|=\overline{\alpha}_1(i)+\cdots+\overline{\alpha}_n(i)$ denotes the total degree of the monomial $e_i$. Then, as explained above, the structure of the matrix \cref{eq:analytic_EDMD} should be close to a lower block-triangular form, that is, $\hat{\mathbf{K}}_{ij} \approx \mathbf{0}$ if $|\overline{\alpha}(i)|<|\overline{\alpha}(j)|$.

A few remarks are in order.
\begin{remark}[Equilibrium]
\label{rem:translated}
The matrix approximation \cref{eq:analytic_EDMD} is based on a Taylor expansion around the origin, and is therefore adapted to an equilibrium at the origin. If the equilibrium $\mathbf{x}^*$ does not lie at the origin, the approximation \cref{eq:analytic_EDMD} can be computed with translated data pairs \mbox{$\{\mathbf{x}_k-\mathbf{x}^*,\mathbf{y}_k-\mathbf{x}^*\}_{k=1}^M$}. This is equivalent to considering a translated kernel \mbox{$k(\mathbf{x}-\mathbf{x}^*,\mathbf{y}-\mathbf{x}^*)$} associated with a RKHS spanned by the orthonormal basis $\{\beta_\alpha(\mathbf{x}-\mathbf{x}^*)^\alpha\}_{\alpha\in \mathbb{N}^n}$. This will be illustrated in \cref{sec:example_1D}. \hfill $\diamond$
\end{remark}
\begin{remark}[Regularization]
    It is well-known that regularization is obtained when the kernel matrix $\mathbf{G}$ is replaced by $\mathbf{G}+\epsilon \mathbf{I}$, where $\epsilon>0$ is the regularization parameter. This is particularly useful in the case of noisy data (see \cref{sec:num_analysis}). \hfill $\diamond$
\end{remark}
\begin{remark}[Non-orthonormal basis]
If the projection subspace $\Hc_N$ is spanned by a non-orthonormal basis of functions $\{e_i\}_{i=1}^N$, the Koopman matrix approximation obtained through orthogonal projection onto $\Hc_N$ is given by
\begin{equation}
\label{eq:approx_non_ortho}
\hat{\mathbf{K}} = (\mathbf{X}^T \mathbf{G}^{-1} \mathbf{X})^{-1} \mathbf{X}^T \mathbf{G}^{-1} \mathbf{Y},
\end{equation}
where $\mathbf{X}^T \mathbf{G}^{-1} \mathbf{X}$ plays the role of a Gram matrix whose entries approximate $\langle e_i,e_j\rangle_\Hc$. The above matrix approximation can be used when the norm of unweighted monomials is not equal to $1$, as an alternative to using weighted monomials. If the basis functions $e_i$ are orthonormal, $\mathbf{e}_i^T \mathbf{G}^{-1} \mathbf{e}_j \approx \langle e_i, e_j \rangle_{\Hc}$ implies that $\mathbf{X}^T \mathbf{G}^{-1} \mathbf{X}\approx \mathbf{I}$ and we recover \cref{eq:analytic_EDMD}. When the basis of monomials is not even orthogonal (e.g. when the kernel is not of Taylor-type), an orthogonal projection will not provide the proper Taylor expansion. Instead, one has to use the dual basis $\{\tilde{e}_i\}_{i=1}^\infty$ such that $\langle e_i,\tilde{e}_j\rangle=\delta_{ij}$. This implies that $\langle \tilde{e}_i, f\rangle = \frac{1}{\overline{\alpha}_1(i)!\cdots \overline{\alpha}_n(i)!}\left.\frac{\partial^{\overline{\alpha}(i)} f}{\partial \mathbf{x}^{\overline{\alpha}(i)}}\right|_{\mathbf{x}=\mathbf{0}}$, and it can be shown that $\tilde{e}_i = \frac{1}{\overline{\alpha}_1(i)!\cdots \overline{\alpha}_n(i)!} \left.\frac{\partial^{\overline{\alpha}(i)} k(\mathbf{x},\cdot)}{\partial \mathbf{x}^{\overline{\alpha}(i)}}\right|_{\mathbf{x}=\mathbf{0}}$  (see e.g. \cite{Ishikawa_JetEDMD}). In this case, the Koopman matrix is given by
\begin{equation}
\label{eq:approx_non_ortho_bis}
\hat{\mathbf{K}} =  \mathbf{\tilde{X}}^T \mathbf{G}^{-1} \mathbf{Y}
\end{equation}
with 
\begin{equation*}
\tilde{\mathbf{X}} = \begin{pmatrix}
\tilde{e}_1(\mathbf{x}_1) & \cdots & \tilde{e}_N(\mathbf{x}_1) \\
\vdots & \ddots & \vdots \\
\tilde{e}_1(\mathbf{x}_M) & \cdots & \tilde{e}_N(\mathbf{x}_M)
\end{pmatrix}
= \begin{pmatrix}
\left.\frac{\partial^{\overline{\alpha}(1)} k(\mathbf{x},\mathbf{x}')}{\partial \mathbf{x}^{\overline{\alpha}(1)}}\right|_{(\mathbf{0},\mathbf{x}_1)} & \cdots & \left.\frac{\partial^{\overline{\alpha}(N)} k(\mathbf{x},\mathbf{x}')}{\partial \mathbf{x}^{\overline{\alpha}(N)}}\right|_{(\mathbf{0},\mathbf{x}_1)} \\
\vdots & \ddots & \vdots \\
\left.\frac{\partial^{\overline{\alpha}(1)} k(\mathbf{x},\mathbf{x}')}{\partial \mathbf{x}^{\overline{\alpha}(1)}}\right|_{(\mathbf{0},\mathbf{x}_M)} & \cdots & \left.\frac{\partial^{\overline{\alpha}(N)} k(\mathbf{x},\mathbf{x}')}{\partial \mathbf{x}^{\overline{\alpha}(N)}}\right|_{(\mathbf{0},\mathbf{x}_M)}
\end{pmatrix}.
\end{equation*}
\hfill $\diamond$
\end{remark}

\paragraph{Related works} We briefly discuss some connections to previous methods. To our knowledge, none of the other data-driven methods exploit the triangular structure of the operator, either due to the choice of projection operator (different from the Taylor projection) or to the choice of subspace (different from a polynomial subspace).

\begin{itemize}
    \item{\textbf{EDMD.}} If the number of data points is equal to the number of basis functions (i.e. $M=N$), the data matrix $\mathbf{X}$ is a square matrix and, provided that it is invertible, it follows from \cref{eq:approx_non_ortho} that $\hat{\mathbf{K}}=\mathbf{X}^{-1} \mathbf{Y}$, which corresponds to standard EDMD with monomials. In this case, the type of projection does not matter and corresponds to polynomial interpolation at the data points (see also \cite{Seenivasaharagavan2025} for other cases). However, when $M>N$, both methods obviously rely on distinct projections.
    \item{\textbf{Kernel-based methods.}} The matrix approximation \cref{eq:analytic_EDMD} is similar to Koopman operator approximations based on kernel ridge regression \cite{Klus_kernel, Kostic}, which also appear in \cite{Bevanda2024} in the context of control theory. In our case, analytic EDMD relies on specific Taylor-type kernels (associated with infinite-dimensional RKHSs). Moreover, contrary to kernel ridge regression, the matrix $\mathbf{X}$ in \cref{eq:analytic_EDMD} is not a (typically infinite-dimensional) feature matrix but should be seen instead as the finite-dimensional basis of the projection subspace. If $\mathbf{X}$ and $\mathbf{Y}$ are interpreted as truncated feature matrices, the kernel trick yields
    \begin{equation}
    \label{eq:kernel_trick}
        \mathbf{YX}^T \approx \mathbf{A}
    \end{equation}
    with
     \begin{equation*}
        \mathbf{A} = \begin{pmatrix}
k(\mathbf{y}_1,\mathbf{x}_1) & \cdots & k(\mathbf{y}_1,\mathbf{x}_M) \\
\vdots & \ddots & \vdots \\
k(\mathbf{y}_M,\mathbf{x}_1) & \cdots & k(\mathbf{y}_M,\mathbf{x}_M)
\end{pmatrix}
\end{equation*}
and it follows that
    \begin{equation}
    \label{eq:ridge_reg}
        \mathbf{X}^T\mathbf{G}^{-1}\mathbf{Yv}=\mu \mathbf{v} \Rightarrow \mathbf{X}^T \mathbf{G}^{-1}\mathbf{YX}^T \mathbf{\tilde{v}}  \approx \mathbf{X}^T \mathbf{G}^{-1} \mathbf{A} \mathbf{\tilde{v}} \approx \mu \mathbf{X}^T \mathbf{\tilde{v}} \Rightarrow  \mathbf{G}^{-1} \mathbf{A} \mathbf{\tilde{v}} \approx \mu  \mathbf{\tilde{v}}
    \end{equation}
    where $\mathbf{v} = \mathbf{X}^T \mathbf{\tilde{v}}$ (see also \cite{Klus_kernel}). This implies that the spectrum of \cref{eq:analytic_EDMD} is equal (in approximation) to the spectrum of the kernel EDMD representation $\mathbf{G}^{-1}\mathbf{A}$ \cite{Williams_kernel}. However, analytic EDMD and kernel EDMD are not equivalent in general since \cref{eq:kernel_trick} and \cref{eq:ridge_reg} are not exact for kernels associated with an infinite-dimensional feature map. These observations are consistent with the fact that kernel EDMD amounts at projecting onto the subspace of reproducing kernel functions $k(\mathbf{x}_i,\cdot)$, and not on a polynomial subspace. Indeed, if we use the (non-orthonormal) basis functions \mbox{$e_i=k(\mathbf{x}_i,\cdot)$}, \cref{eq:data_mat} and \cref{eq:approx_non_ortho} yield $\mathbf{X=G}$, $\mathbf{Y=A}$, and $\hat{\mathbf{K}}=\mathbf{G}^{-1} \mathbf{A}$, respectively, so that we recover kernel EDMD.
    \item{\textbf{JetEDMD.}} Analytic EDMD bears similarity to the recent JetEDMD \cite{Ishikawa2024, Ishikawa_JetEDMD}, which also uses a truncation method with a basis of monomials (i.e. Taylor projection). However, JetEDMD still makes use of standard EDMD in its first step, and therefore relies on both Galerkin and Taylor projections, while analytic EDMD solely relies on the Taylor projection.
\end{itemize}

\subsection{Spectral properties}

The spectral properties of the Koopman operator can be approximated by the eigenvalues and eigenvectors of the Koopman matrix $\cref{eq:analytic_EDMD}$. In particular, the eigenvalues of $\hat{\mathbf{K}}$ provide an approximation of the eigenvalues of $K$ and the associated eigenvectors yield the Taylor coefficients of the corresponding eigenfunctions. Note that the approximation error on the eigenvalues and Taylor coefficients of the eigenfunctions is only due to the data-driven setting, and not to the finite-dimensional nature of the method.

The eigenvalues and eigenvectors of the Koopman matrix can be computed by exploiting the block-triangular structure, following similar lines as in \cite{MauroyMezic_CDC}. Let us denote by $\overline{\mathbf{K}}_{rs}$, with \mbox{$r,s\in\{0,\dots,|\alpha|_{max}\triangleq\overline{\alpha}(N)\}$}, the matrix block of $\hat{\mathbf{K}}$ containing all entries $K_{ij}$ such that $\overline{\alpha}(i)=r$ and $\overline{\alpha}(j)=r$.
\begin{itemize}
\item{\textbf{Eigenvalues.}} The eigenvalues of $\hat{\mathbf{K}}$ are obtained by computing the eigenvalues $\hat{\mu}$ of the diagonal matrix blocks $\overline{\mathbf{K}}_{rr}$. They are the estimates of the Koopman eigenvalues $\mu=\mu_1^{\alpha_1} \cdots \mu_n^{\alpha_n}$, with $|\alpha|=r$ and where $\mu_j$ are the eigenvalues of the Jacobian matrix $\mathbf{J}$ of the map $\varphi$ evaluated at the equilibrium \cite{Maccluer1984spectra}. Note that the eigenvalues $\hat{\mu}_j$ of $\overline{\mathbf{K}}_{11} \approx \mathbf{J}^T$ are the estimates of $\mu_j$.
\begin{remark}[Continuous time]
    In the continuous-time setting, the eigenvalues $\lambda$ of the Koopman generator are estimated by computing $\hat{\lambda}=\log(\hat{\mu})/\Delta t$ where $\Delta t$ is the sampling time, with $\Delta t<\pi/\Im\{\lambda\}$ to prevent aliasing \cite{Zeng_sampling}. In this case, the lattice-structured spectrum of the Koopman generator consists of eigenvalues of the form \mbox{$\lambda=\alpha_1 \lambda_1 + \cdots + \alpha_n \lambda_n$}, where $\lambda_j$ are the eigenvalues of the Jacobian matrix $\mathbf{J}$ of the vector field evaluated at the equilibrium (see e.g. \cite{Mezic_ann_rev}). \hfill $\diamond$
\end{remark}
\item{\textbf{Eigenfunctions.}} We focus on the eigenfunctions associated with the eigenvalues $\mu_j \in \sigma(\mathbf{J})$, i.e. the case $|\alpha|=1$. These are the principal eigenfunctions which completely capture the dynamics. Other eigenfunctions are redundant and correspond to integer powers of the principal ones. Let $\overline{\mathbf{v}}_r^{(j)}$ be the components of the eigenvector of $\hat{\mathbf{K}}$ associated with $\hat{\mu}_j$, related to monomials of total degree $r$, i.e. $\overline{\mathbf{v}}_r^{(j)}$ contains the $r$-th order approximate Taylor coefficients (in the basis $\{e_i\}_{i=1}^N$) of $\phi_{\mu_j}$. It follows from the block-triangular structure of $\hat{\mathbf{K}}$ that
\begin{equation*}
\sum_{s=1}^{r} \overline{\mathbf{K}}_{rs} \overline{\mathbf{v}}_s^{(j)} = \hat{\mu}_j \overline{\mathbf{v}}_r^{(j)}.
\end{equation*}
We can then proceed recursively. For $r=1$, $\overline{\mathbf{v}}_1^{(j)}$ is the eigenvector of $\overline{\mathbf{K}}_{11} \approx \mathbf{J}^T$ associated with $\hat{\mu}_j$ and, for $r>1$, we have
\begin{equation}
\label{eq:eigvec}
\overline{\mathbf{v}}_r^{(j)} = (\overline{\mathbf{K}}_{rr} - \hat{\mu}_j \mathbf{I})^{-1} \sum_{s=1}^{r-1} \overline{\mathbf{K}}_{rs} \overline{\mathbf{v}}_s^{(j)}.
\end{equation}
Once the Taylor coefficients are computed, an approximation of the eigenfunction is given by
\begin{equation}
\label{eq:taylor_Koop_eigfct}
\hat{\phi}_{\mu_j} = \sum_{r=1}^{|\alpha|_{max}} (\overline{\mathbf{v}}_r^{(j)})^T \, \mathbf{\overline{e}}_r,
\end{equation}
where $\mathbf{\overline{e}}_r$ is the vector of monomials $e_i$ of total degree equal to $r$. It should be noted that the validity of the approximation of the eigenfunction is limited by the radius of convergence of the Taylor series, which depends on the location of the (possibly complex-valued) equilibria of the dynamics (see the discussion in \cite{Mauroy_Mezic_stability}).
\end{itemize} 

\subsection{Algorithm}

Analytic EDMD is summarized in \cref{alg:analytic_EDMD}. We remark that the accuracy of the Koopman matrix can be verified with the block-triangular structure. Moreover, one can check whether the eigenvalues of the matrix match the known lattice structure of the exact eigenvalues, a method that can potentially be useful to estimate the maximal order for which the data-driven Taylor approximation remains accurate. Note that error bounds will be given in \cref{sec:errors}.

\begin{algorithm}[h]
	\caption{Analytic EDMD}
	\label{alg:analytic_EDMD}
	\begin{algorithmic}[1]
\Statex{\bf Input:} Snapshot pairs $\{(\mathbf{x}_k,\mathbf{y}_k)\}_{k=1}^M$ (possibly with sampling time $\Delta t$); equilibrium point $\mathbf{x}^*$; Taylor-type kernel function $k$; basis of (weighted) monomials $\{e_i\}_{i=1}^N$ (with maximal total degree $|\alpha|_{max}$); regularization parameter $\epsilon$.
\Statex{\bf Output:} Koopman matrix approximation $\hat{\mathbf{K}}$; set $S_{d}$ or $S_{c}$ of (discrete or continuous-time) Koopman eigenvalues; principal Koopman eigenfunctions $\hat{\phi}_{\mu_j}$
\State $\mathbf{x}_k \leftarrow \mathbf{x}_k-\mathbf{x}^*$, $\mathbf{y}_k \leftarrow \mathbf{y}_k-\mathbf{x}^*$ ($k=1,\dots,M$)
\newline
\textbf{Koopman matrix approximation}
\State Construct the $M \times M$ kernel matrix $\mathbf{G}$ defined in \cref{eq:kernel_mat}
\State Construct the $M\times N$ data matrices $\mathbf{X}$ and $\mathbf{Y}$ defined in \cref{eq:data_mat}
\State Compute the $N \times N$ matrix $\hat{\mathbf{K}} = \mathbf{X}^T (\mathbf{G+\epsilon \mathbf{I}})^{-1} \mathbf{Y}$ \newline (or $\hat{\mathbf{K}} = (\mathbf{X}^T (\mathbf{G+\epsilon \mathbf{I}})^{-1} \mathbf{X})^{-1} \mathbf{X}^T (\mathbf{G+\epsilon \mathbf{I}})^{-1} \mathbf{Y}$ if the basis $\{e_i\}$ is not orthonormal)
\newline
\textbf{Koopman eigenvalues}
\For {$r=1:|\alpha|_{max}$} 
			\State Compute the eigenvalues of the matrix block $\overline{\mathbf{K}}_{rr}$ and add them to the set $S_d$
\EndFor
\State For continuous-time systems, add $\hat{\lambda}=\log(\hat{\mu})/\Delta t$ to the set $S_c$ for all $\hat{\mu}\in S_d$
\newline
\textbf{Principal Koopman eigenfunctions}
\State Compute the eigenvalues $\hat{\mu}_j$ and eigenvectors $\mathbf{w}_j$ of the matrix block $\overline{\mathbf{K}}_{11}$
\For {$j=1:n$} 
\State $\overline{\mathbf{v}}_1^{(j)} \leftarrow \mathbf{w}_j$
\For {$r=2:|\alpha|_{max}$} 
			\State $\overline{\mathbf{v}}_r^{(j)} \leftarrow (\overline{\mathbf{K}}_{rr} - \hat{\mu}_j \mathbf{I})^{-1} \sum_{s=1}^{r-1} \overline{\mathbf{K}}_{rs} \overline{\mathbf{v}}_s^{(j)}$
\EndFor
\EndFor
\State Compute the principal Koopman eigenfunctions $\hat{\phi}_{\mu_j}$ using \cref{eq:taylor_Koop_eigfct}.
	\end{algorithmic}
\end{algorithm}

\section{Numerical simulations}
\label{sec:num_ex}

Analytic EDMD is illustrated with several numerical simulations in the case of continuous-time systems that possess one (or several) hyperbolic equilibrium. We will also consider specific situations where the data points are taken far from the equilibrium, and only partial state measurement is allowed. In all cases, we will use the Szeg\"o kernel \mbox{$k(\mathbf{x},\mathbf{y})=\prod_{i=1}^n 1/(1-x_i y_i)$} with no regularization ($\epsilon=0$).

\subsection{First examples}

\begin{figure}
\centering
\subfigure[Unstable eigenvalues (at $x^*=0$)]{\includegraphics[width=0.45\linewidth]{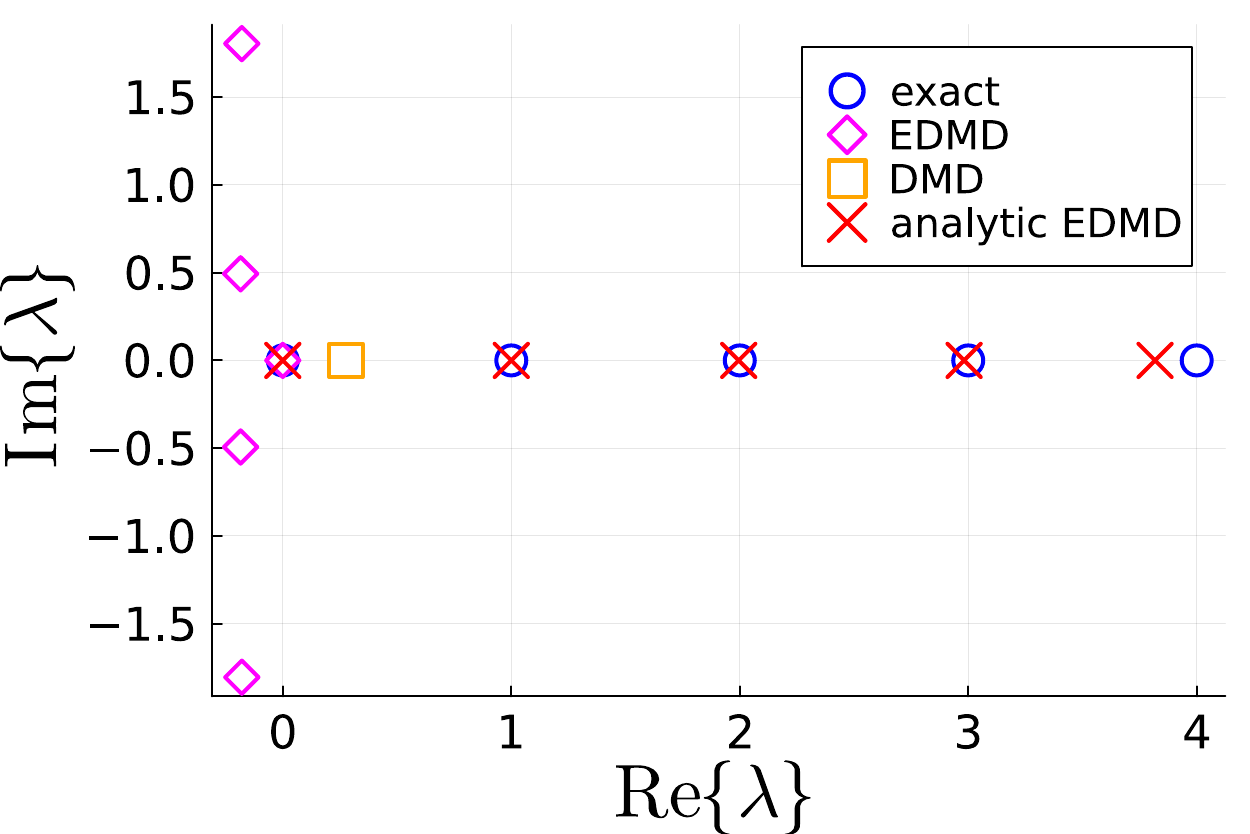}}
\subfigure[Stable eigenvalues (at $x^*=1$)]{\includegraphics[width=0.45\linewidth]{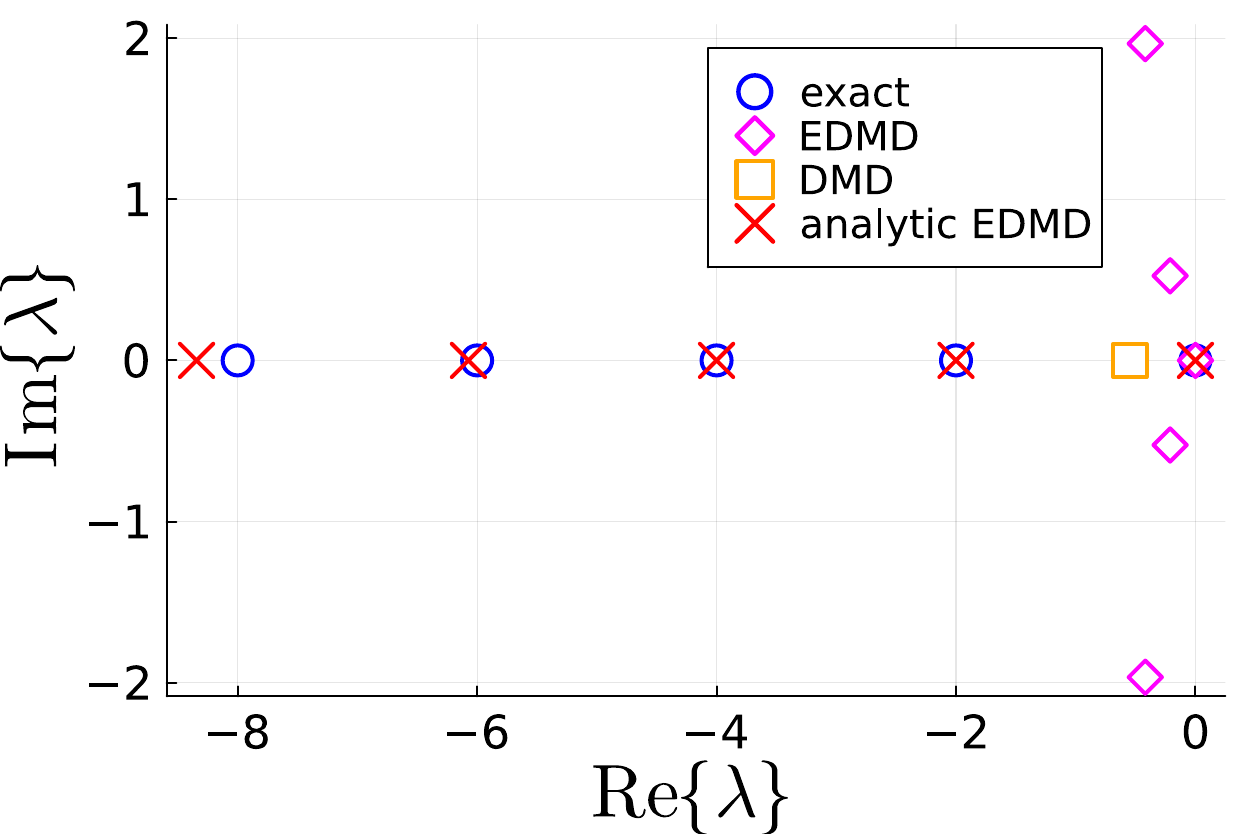}}
\caption{Computation of Koopman eigenvalues for a one-dimensional cubic dynamics.}
\label{fig_1D}
\end{figure}

\begin{example}
\label{sec:example_1D}
Let us consider the one-dimensional cubic dynamics $\dot{x} = x-x^3$, which possess one unstable equilibrium at $x^*=0$ and two stable equilibria at $x^*=-1$ and \mbox{$x^*=1$}. The Koopman eigenvalues are computed with DMD, EDMD, and analytic EDMD with monomial basis functions up to degree $|\alpha|_{max}=4$. A set of $M=20$ data pairs is generated from a uniform distribution over the interval $X=[0,1]$, with the sampling time $\Delta t = 0.5$. The results are shown in \cref{fig_1D}(a). We note that unstable eigenvalues $\lambda_j=j$, $j\in \mathbb{N}$ are computed with analytic EDMD since the projection on monomials is related to the Taylor expansion around the unstable equilibrium $x^*=0$. By using the translated kernel $1/(1-(x-1)(y-1))$ associated with a Taylor expansion around the stable equilibrium at $x^*=1$ (see \cref{rem:translated}), analytic EDMD can also capture stable Koopman eigenvalues $\lambda_j=-2j$, $j\in \mathbb{N}$, with the same data set (see \cref{fig_1D}(b)). As mentioned in \cite{Page_Kerswell}, the performance of EDMD critically depends on the region where the data points are sampled. In the present case of a large sampling region, EDMD captures neither stable nor unstable Koopman eigenvalues while analytic EDMD demonstrates good performance.
\end{example}

\begin{figure}[h!]
\centering
\subfigure[$M=50$, $\Delta t=1$]{\includegraphics[width=0.45\linewidth]{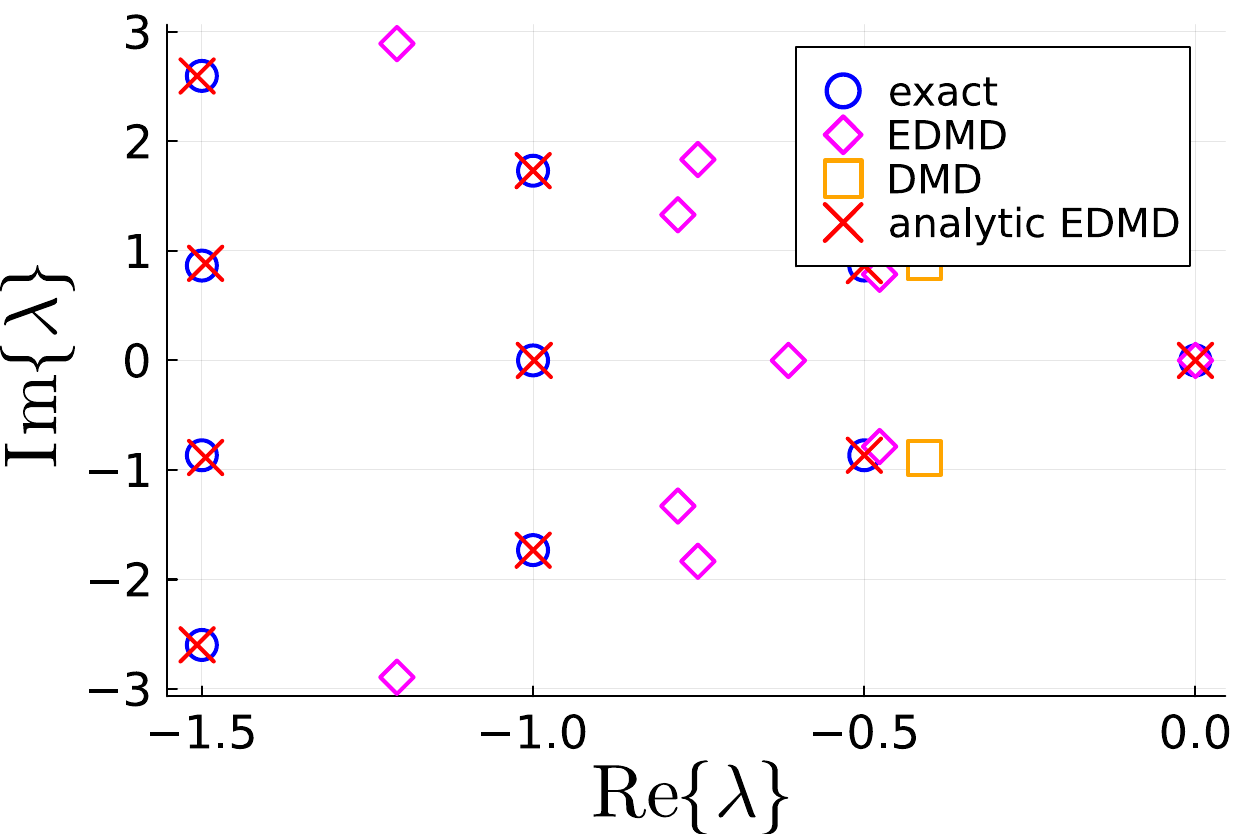}}
\subfigure[$M=250$, $\Delta t=0.5$]{\includegraphics[width=0.45\linewidth]{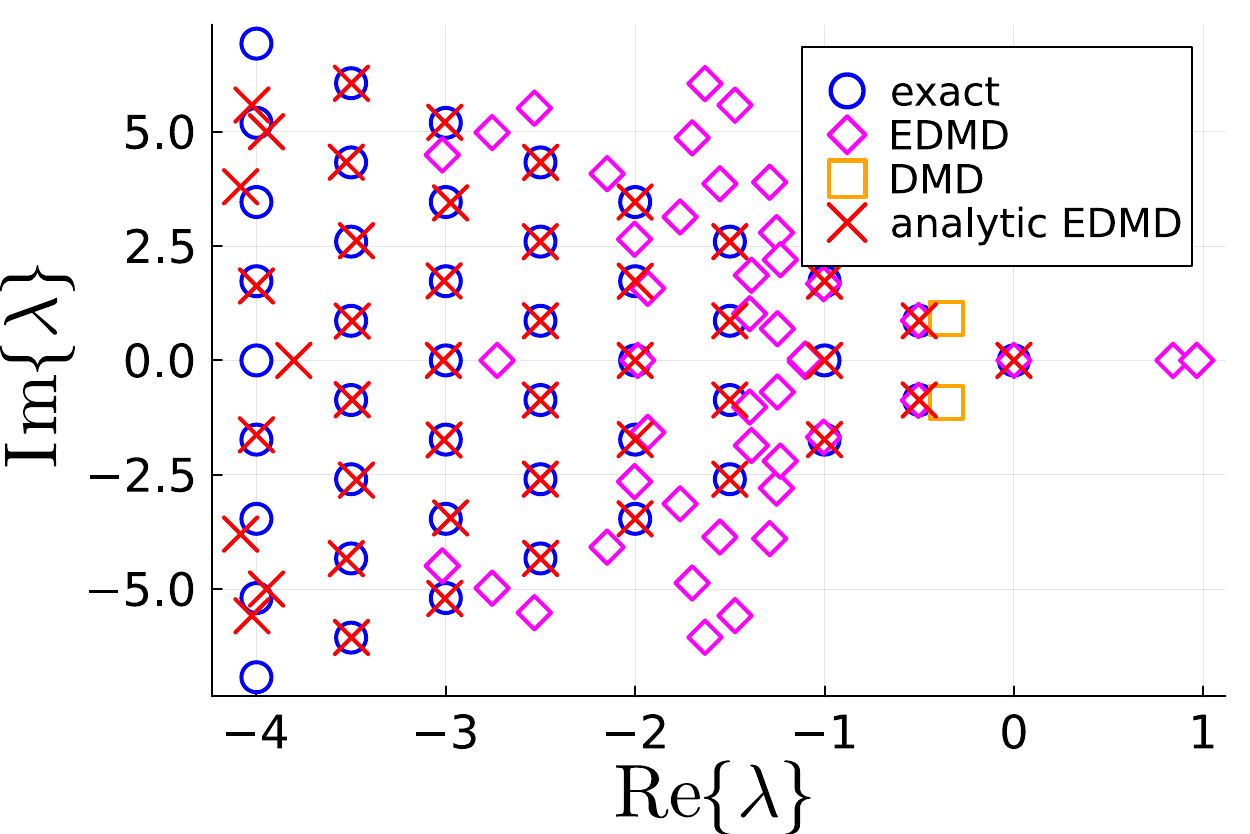}}
\caption{Computation of Koopman eigenvalues for the Van der Pol dynamics.}
\label{fig_VDP}
\end{figure}

\begin{figure}
\centering
\subfigure[Absolute value]{\includegraphics[width=0.45\linewidth]{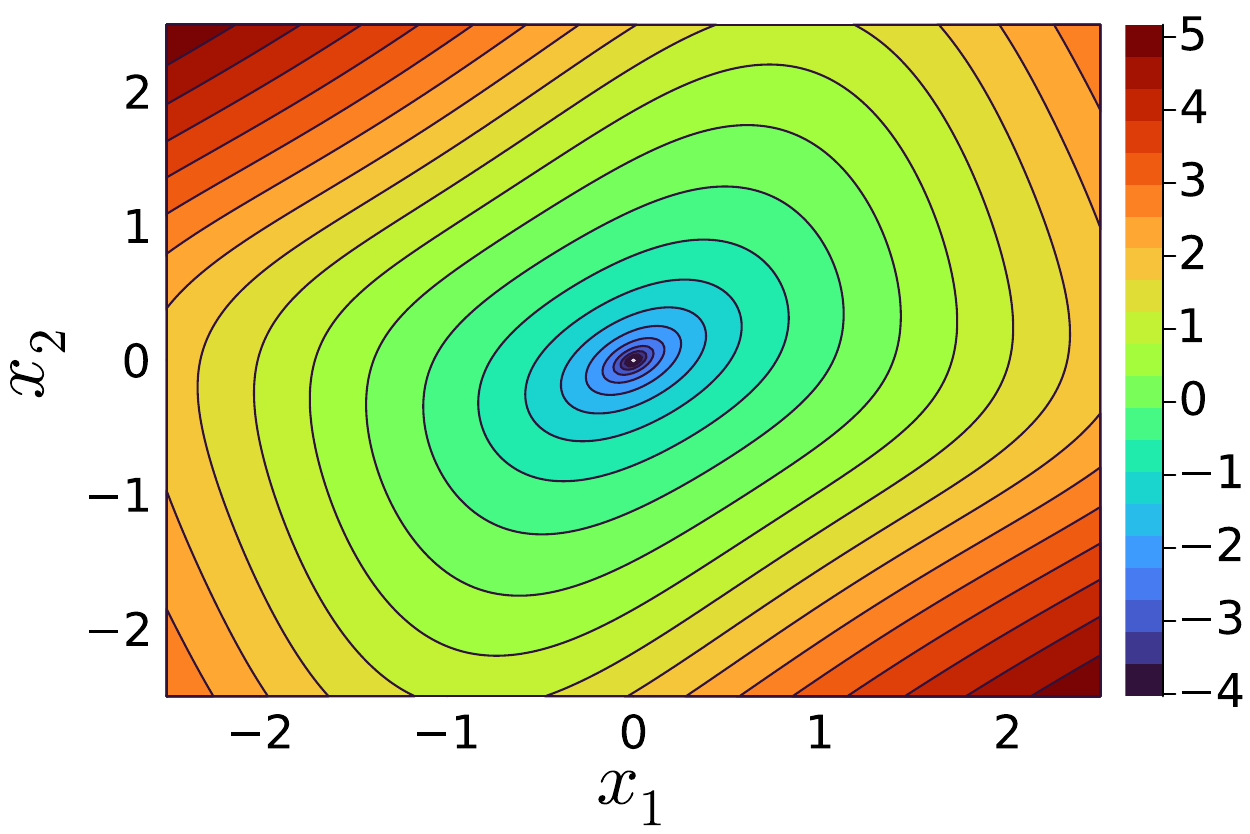}}
\subfigure[Argument]{\includegraphics[width=0.45\linewidth]{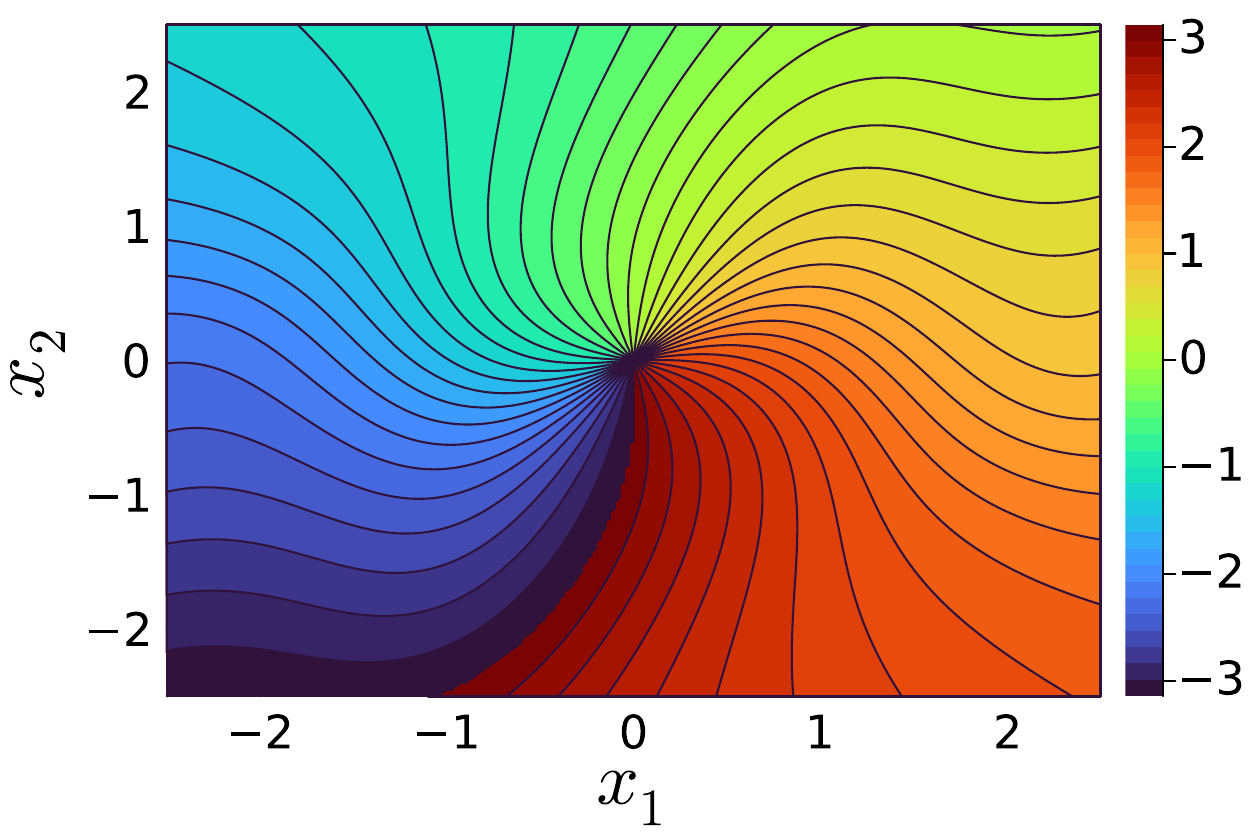}}
\caption{Taylor approximation (to the 8th order) of a principal Koopman eigenfunction for the stable Van der Pol dynamics. (The colormap shows the logarithm of the absolute value (left) and the value of the argument (right).)}
\label{fig_VDP_eigfct}
\end{figure}

\begin{example}
\label{sec:Van_der_Pol}

We consider the Van der Pol dynamics 
\begin{eqnarray*}
\dot{x}_1 & = & -x_2 \\
\dot{x}_2 & = & -(1-x_1^2)x_2+x_1
\end{eqnarray*}
which possess a stable equilibrium at the origin and an unstable limit cycle. The eigenvalues of the Jacobian matrix of the system at the origin are $-0.5\pm i\sqrt{3}/2$. The Koopman eigenvalues are computed with DMD, EDMD, and analytic EDMD in two different settings. In the first setting, $M=50$ data pairs are uniformly distributed over $X=[-1,1]^2$ with the sampling time $\Delta t=1$, and monomial basis functions are used up to total degree $|\alpha|_{max}=3$ (\cref{fig_VDP}(a)). In the second setting, $M=250$ data pairs are generated over $[-1,1]^2$ with the sampling time $\Delta t=0.5$, and monomial basis functions are used up to total degree $|\alpha|_{max}=8$ (\cref{fig_VDP}(b)). In both cases, we observe that analytic EDMD yields very accurate results. In particular, there is no spurious eigenvalue, even with a small dataset (\cref{fig_VDP}(a)) or when fast eigenvalues are captured (\cref{fig_VDP}(b)). The principal eigenfunction associated with the eigenvalue $\lambda=-0.5+ i\sqrt{3}/2$ is also computed in the second setting. The level sets of the absolute value of the eigenfunction  are the so-called isostables of the equilibrium \cite{MMM_isostables} (\cref{fig_VDP_eigfct}(a)), while the level sets of the argument correspond to a periodic partition of the basin of attraction (\cref{fig_VDP_eigfct}(b)). The results are consistent with those obtained in \cite{MauroyMezic_CDC} from the vector field.
\end{example}

\begin{figure}[h!]
\centering
\subfigure[First principal eigenfunction]{\includegraphics[width=0.45\linewidth]{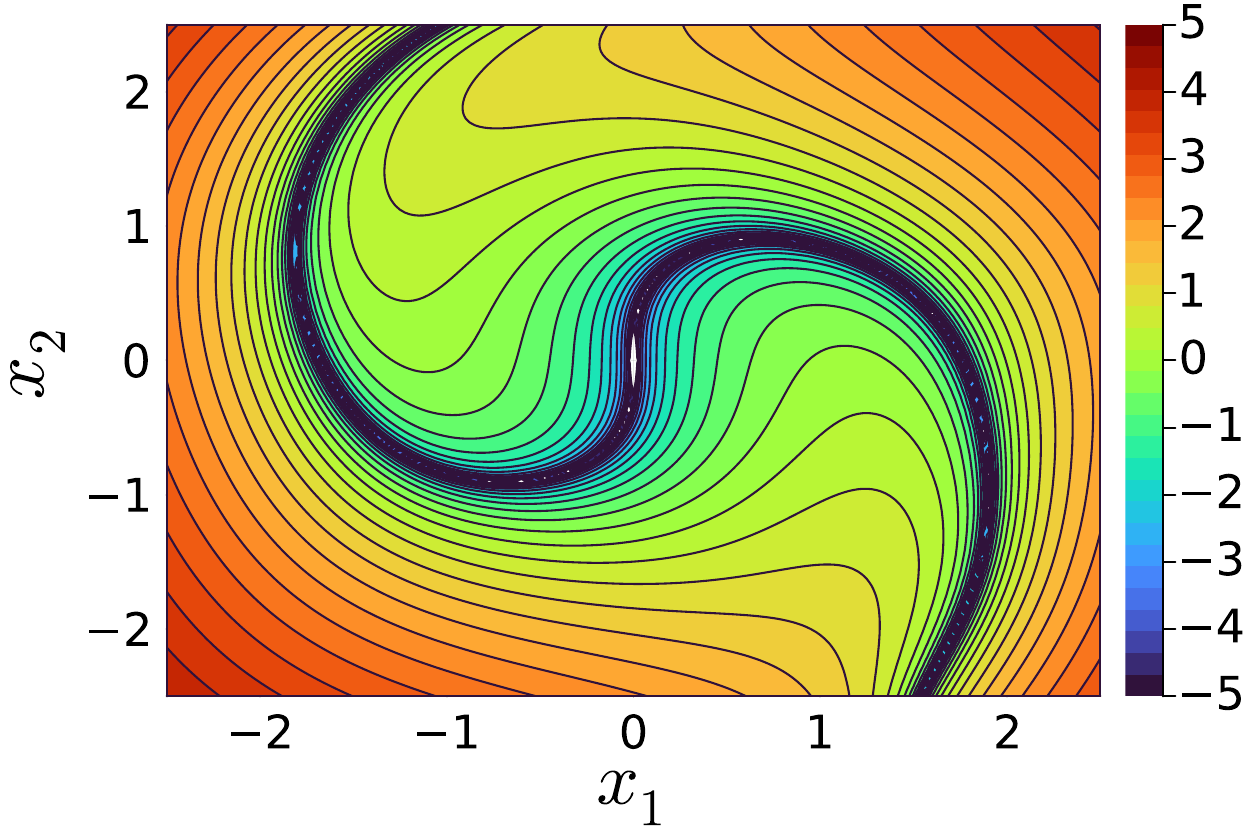}}
\subfigure[Second principal eigenfunction]{\includegraphics[width=0.45\linewidth]{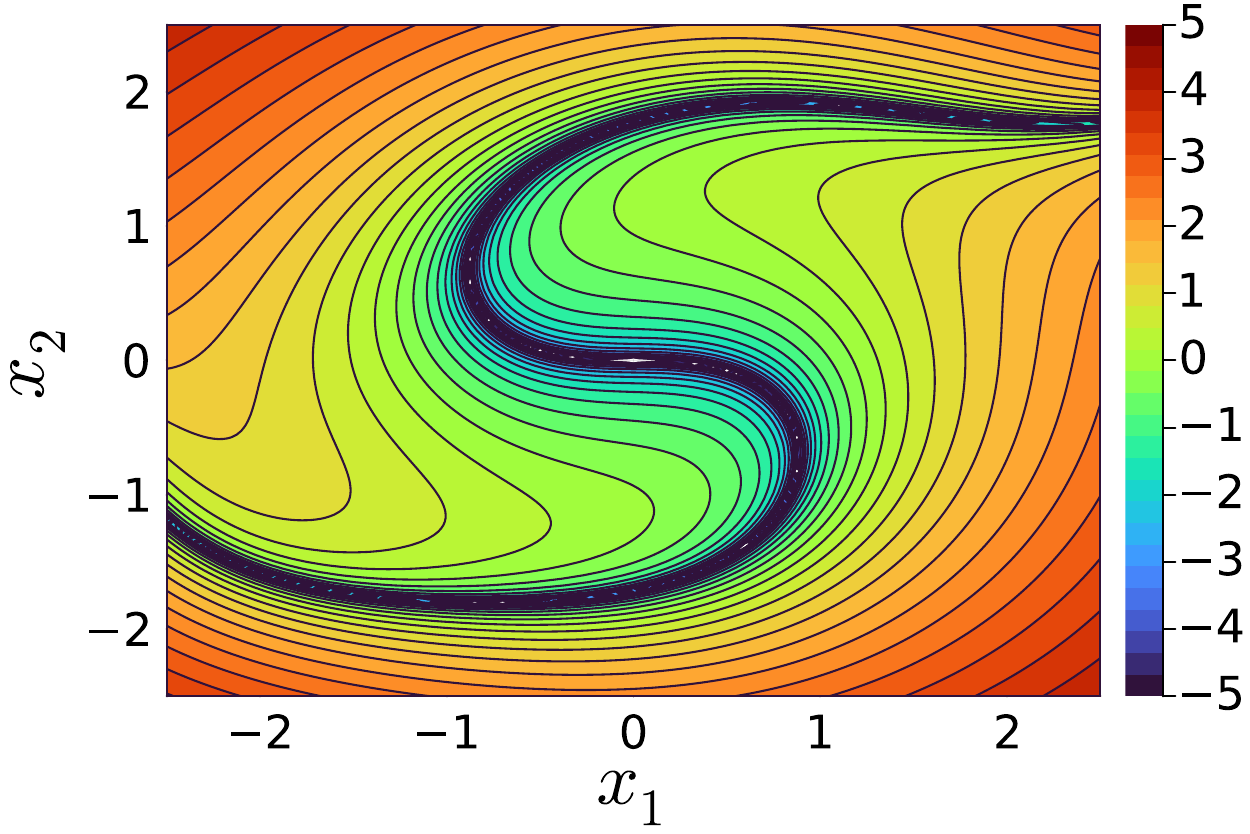}}
\caption{Taylor approximation (to the 6th order) of the principal Koopman eigenfunctions for the nonlinear damped oscillator. (The colormap shows the logarithm of the absolute value of the eigenfunctions.)}
\label{fig_rot_eigfct}
\end{figure}

\begin{example}
The nonlinear damped oscillator $\dot{r} = -r$, $\dot{\theta}=r^2$ in polar coordinates can be rewritten as
\begin{eqnarray*}
\dot{x}_1 & = & -x_1-x_1^2 x_2-x_2^3 \\
\dot{x}_2 & = & -x_2+x_1x_2^2+x_1^3
\end{eqnarray*}
in Cartesian coordinates. The dynamics possess a globally stable equilibrium at the origin and the linearized system admits a single eigenvalue $\lambda=-1$. We generate $M=50$ data pairs, which are uniformly randomly distributed on $X=[-1,1]^2$, with the sampling time $\Delta t=2$. Two principal Koopman eigenfunctions are computed with analytic EDMD (with monomial basis functions up to total degree $6$). The results are shown in \cref{fig_rot_eigfct} and are in agreement with the results obtained in \cite{MauroyMezic_CDC} from the vector field.
\end{example}

\subsection{Specific cases}
\label{subsec:specific_cases}

We now show that analytic EDMD is efficient when data points are taken far from the equilibrium and when only partial state measurements are allowed.\\

\paragraph*{Data points far from the equilibrium}

\begin{figure}[h!]
\centering
\subfigure[Data points]{\includegraphics[width=0.45\linewidth]{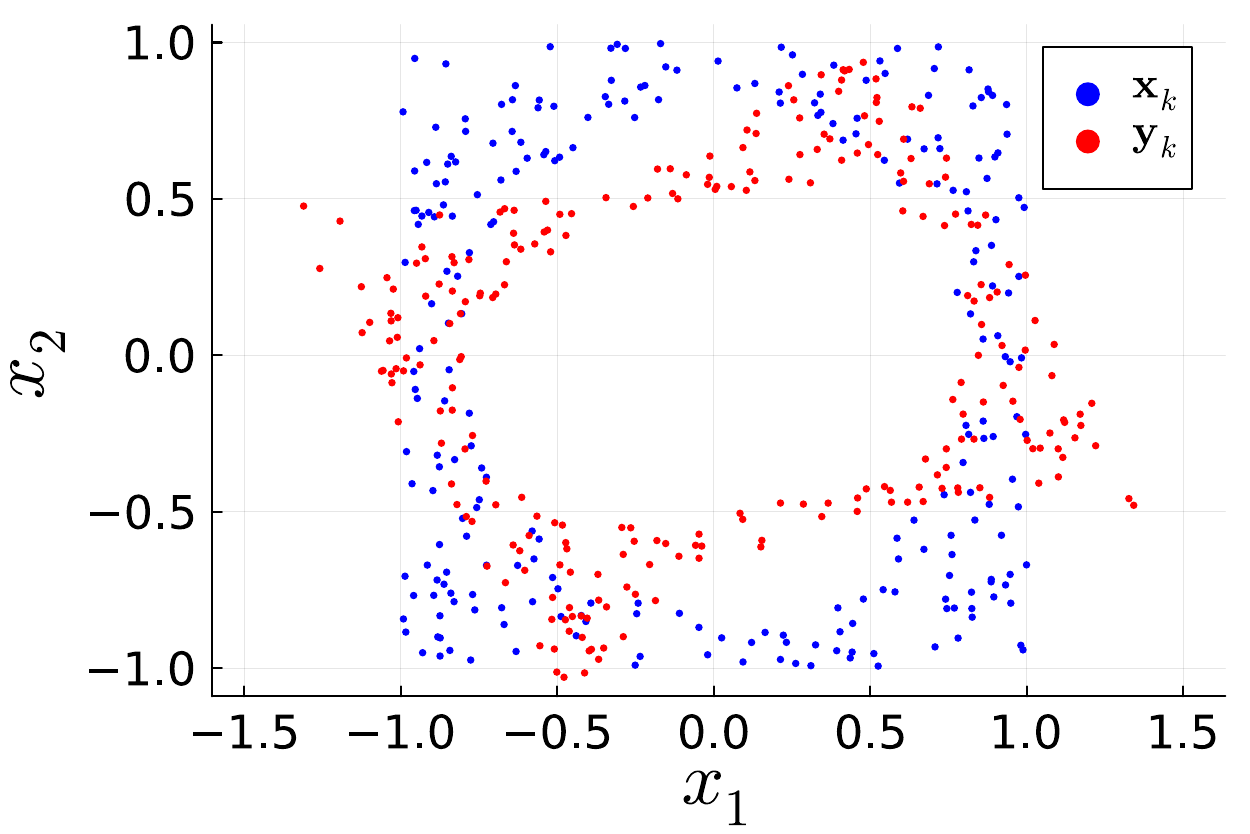}}
\subfigure[Koopman eigenvalues]{\includegraphics[width=0.45\linewidth]{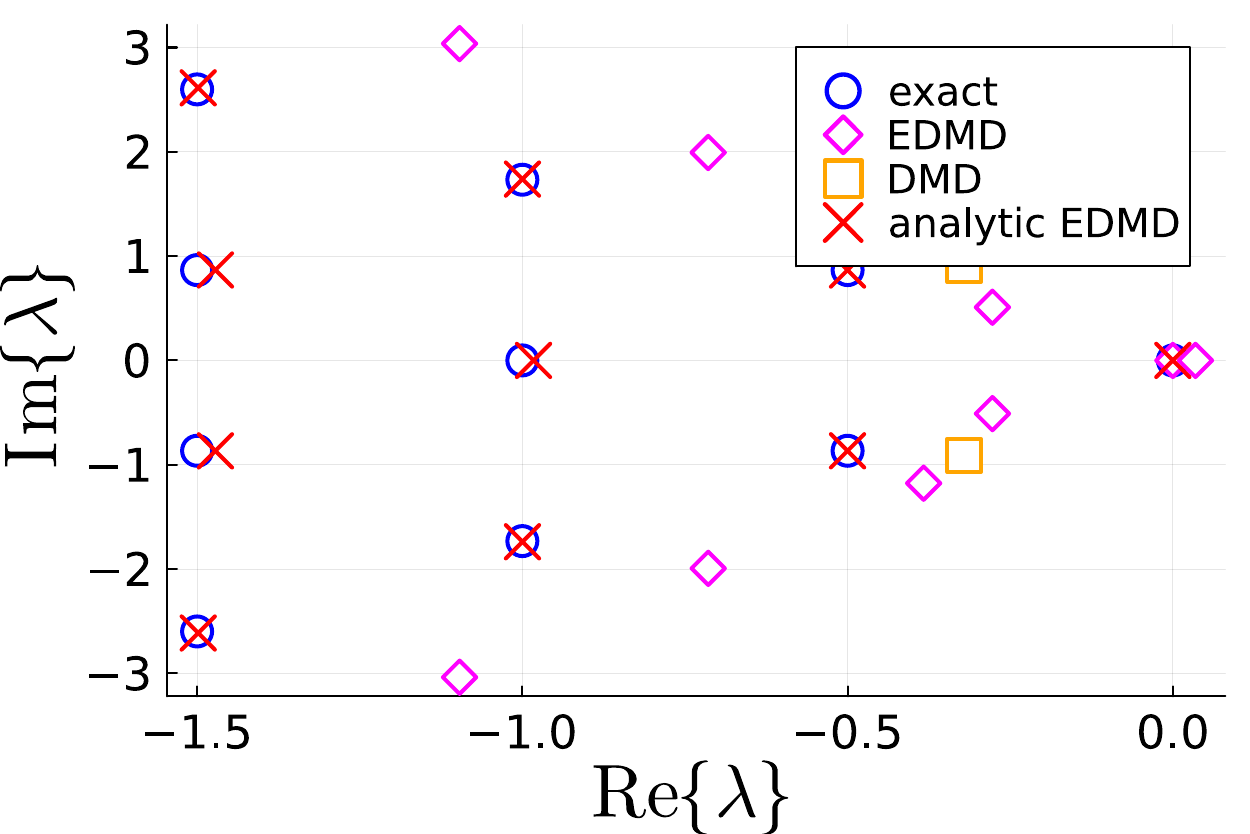}}
\caption{Analytic EDMD is efficient to compute the Koopman eigenvalues, even when data points are taken far from the equilibrium. The Van der Pol dynamics are considered here.}
\label{fig_far_data_points}
\end{figure}

We consider again the Van der Pol dynamics, in a similar setting as in \cref{sec:Van_der_Pol} ($M=250$ data pairs, sampling time $\Delta t=0.5$, monomial degree up to $|\alpha|_{max}=3$), but the data points $\mathbf{x}_k$ are generated over $X=[-1,1]^2$ with the constraint $\|\mathbf{x}_k\|> 0.8$ (\cref {fig_far_data_points}(a)). In this case, it is remarkable that analytic EDMD is still very efficient to capture the spectral properties of the Koopman operator, including the eigenvalues of the Jacobian matrix at the equilibrium (\cref{fig_far_data_points}(b)). Such performance can be explained by the fact that the results are based on the approximation of real-analytic functions, which are uniquely defined by their values on an arbitrary set of positive Lebesgue measure, according to the identity theorem (see e.g. \cite{identity_theorem}). This property is confirmed in \cref{prop:bound_entries} below. \\

\begin{figure}[h!]
\centering
\includegraphics[width=0.5\linewidth]{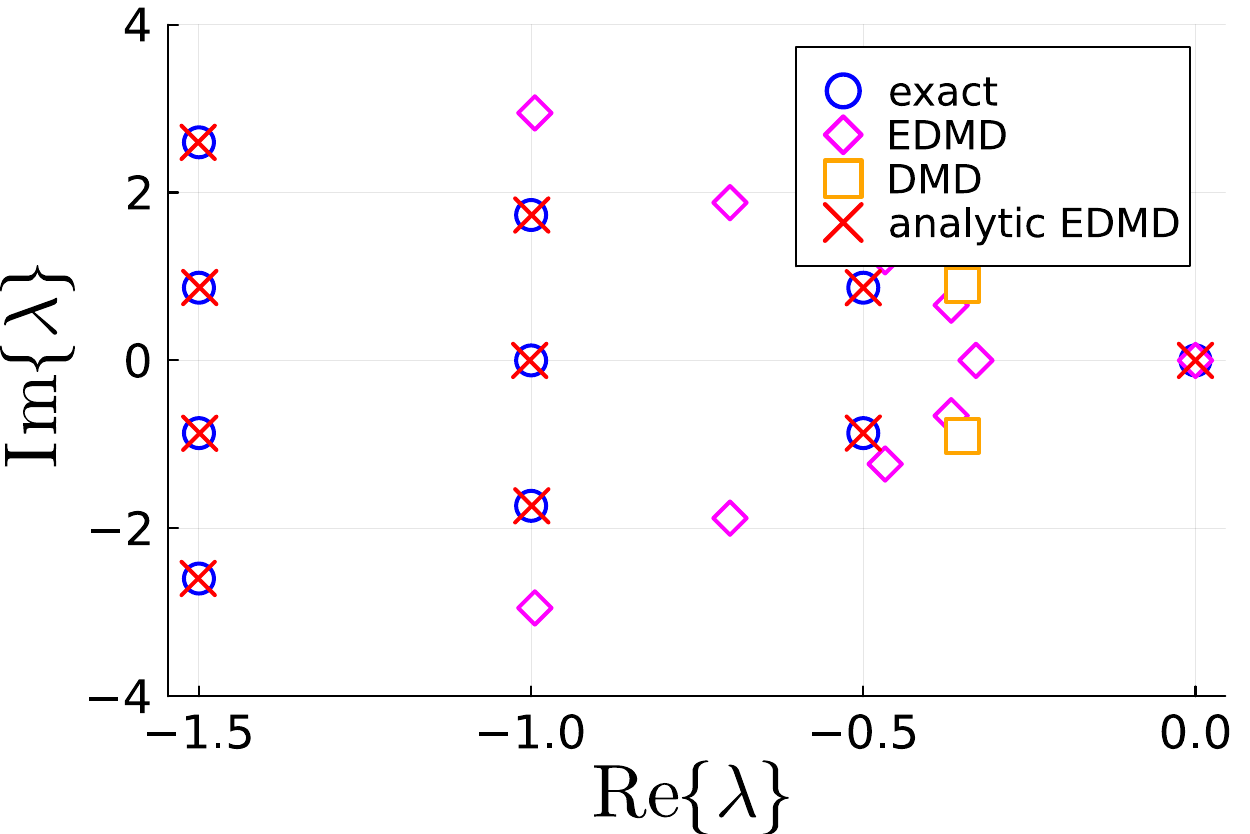}
\caption{Computation of Koopman eigenvalues for the Van der Pol dynamics with time delays, when one state is measured.}
\label{fig_VDP1_partial}
\end{figure}

\paragraph*{Partial state measurement}
Now, we consider the case of partial state measurements and replace non-measured states by time delays of measured states (delay embedding). According to Takens' embedding theorem, the time-delay dynamics are described by a map $\varphi_d$ which satisfies $\varphi_d \circ \Phi = \Phi \circ \varphi_d$, where $\Phi$ is a delay embedding, and therefore inherit the spectral properties of the original dynamics. This motivates the use of Prony-based techniques \cite{Yoshi_prony} and Hankel-EDMD methods \cite{Arbabi} to capture the spectral properties of the Koopman operator (see also \cite{Das_delay, Brunton_time_delays, Koltai_Koopman_Takens} for connections between the Koopman operator and delay embeddings). Analytic EDMD is illustrated in this context with the Van der Pol dynamics, where only the first state is measured. The sampled trajectories $\{\mathbf{x}_k,\varphi(\mathbf{x}_k),\varphi^{2}(\mathbf{x}_k)\}$ are generated for $100$ initial conditions randomly distributed on $X=[0,1]^2$ with the sampling time $\Delta t=0.5$, resulting in $M=100$ snapshot data pairs of time-delayed coordinates $\{([\mathbf{x}_k]_1,[\varphi(\mathbf{x}_k)]_1), ([\varphi(\mathbf{x}_k)]_1,[\varphi^{2}(\mathbf{x}_k)]_1)\}_{k=1}^M$ where $[\mathbf{x}_k]_1$ denotes the first component of $\mathbf{x}_k$. Monomial basis functions are used up to total degree $|\alpha|_{max}=3$. The results depicted in \cref{fig_VDP1_partial} show that the Koopman spectrum can be efficiently recovered with partial measurements. 

\section{Numerical analysis}
\label{sec:num_analysis}

In this section, we perform a few numerical tests with continuous-time systems to assess the performance of analytic EDMD. We define several performance metrics as follows. Consider the set $\sigma_r=\{ \Pi_{j=1}^n \lambda_j^{\alpha_j} : |\alpha|=r\}$ of exact (continuous-time) eigenvalues associated with the order $r$, the exact spectrum $\sigma=\cup_{k\geq 0} \sigma_k$, and the set $S_c$ of (continuous-time) eigenvalues computed with a given numerical method. We define the estimated spectrum accuracy (ESA) as the greatest distance between any exact eigenvalue (associated with the order $r$) and the set of estimated eigenvalues, that is
\begin{equation}
\label{eq:metric1}
    \mathrm{ESA}_r = \max_{\lambda \in \sigma_r} \min_{\hat{\lambda} \in S_c} |\hat{\lambda}-\lambda|.
\end{equation}
In contrast, the spectral pollution measure (SPM) identifies spurious eigenvalues through the average closest distance between any estimated eigenvalue and the exact spectrum, that is
\begin{equation}
\label{eq:metric2}
    \mathrm{SPM} = \frac{1}{\# S_c} \sum_{\hat{\lambda} \in S_c} \min_{\lambda \in \sigma} |\hat{\lambda}-\lambda|,
\end{equation}
where $\# S_c$ denotes the cardinality of $S_c$. Finally, the accuracy of a computed eigenfunction $\hat{\phi}_{\lambda}$ will be measured through the ratio $\hat{\phi}_\lambda(\varphi(\mathbf{x})) / \hat{\phi}_\lambda(\mathbf{x}) \approx e^{\lambda \Delta t}$, where $\varphi$ is the flow map generated by the dynamics and evaluated at time $\Delta t$. In particular, we define the estimated eigenfunction accuracy (EFA) as the mean relative error
\begin{equation}
\label{eq:metric3}
   \mathrm{EFA}  = \frac{1}{M_{test}} \sum_{k=1}^{M_{test}}  \frac{\left|\frac{\hat{\phi}_\lambda(\varphi(\mathbf{x}_k))}{\hat{\phi}_\lambda(\mathbf{x}_k)}-e^{\lambda \Delta t} \right|}{\left|e^{\lambda \Delta t}\right|}
\end{equation}
where $\{\mathbf{x}_k\}_{k=1}^{M_{test}}$ are test points uniformly randomly distributed over $X$. We will choose $M_{test}=50$ and consider the estimated principal eigenfunction associated with the closest estimate $\hat{\lambda}$ of the dominant Jacobian eigenvalue $\lambda_1$ (i.e. $\Re\{\lambda_1\} = \max_{\lambda_j \in \sigma(\mathbf{J})} \Re\{\lambda_j\}$).

We will consider several test cases: 
\begin{enumerate}
    \item the stable Van der Pol dynamics (see \cref{sec:Van_der_Pol});
    \item the Duffing dynamics (with a single trajectory);
\begin{eqnarray*}
\dot{x}_1 & = & x_2 \\
\dot{x}_2 & = & -0.5\,x_2-x_1(-1+x1^2)
\end{eqnarray*}
which possesses two stable hyperbolic equilibria at $(-1,0)$ and $(1,0)$,
\item a $10$-dimension stable quadratic networked dynamics of the form
\begin{equation*}
    \dot{\mathbf{x}} = \mathbf{J} \mathbf{x} -0.2 \, \mathbf{x} \odot \mathbf{x},\qquad \mathbf{x}\in \mathbb{R}^{10},
\end{equation*}
where $\mathbf{J}$ is a Hurwitz matrix whose off-diagonal and diagonal entries are uniformly randomly distributed over $[-1,0]$ and $[-2,-1]$, respectively, and where $\odot$ denotes the Hadamard (component-wise) product. 
\end{enumerate}
The simulation parameters are summarized in \cref{tab:simu_param}.

\begin{table}[h!]
    \centering
    \begin{tabular}{cccccc}
        Dynamics & $M$ & $X$ & $\Delta t$ & $|\alpha|_{max}$ & data\\
        \hline
       Van der Pol & $75/250$ & $[-1,1]^2$ & $0.5$ & $6$ & $M$ trajectories \\
       Duffing  & $100/250$ & $[-1,1]^2$ & $0.1$ & $3$ & one trajectory \\
       $10$-dim. network & $1100$ & $[-0.3,0.3]^2$ & $0.5$ & $2$ & $M$ trajectories \\
    \end{tabular}
    \caption{Simulation parameters}
    \label{tab:simu_param}
\end{table}

\subsection{Comparison with other methods}
\label{sec:comparison}

We compare the performance of analytic EDMD with related methods mentioned in \cref{sec:finite_dim_approx} (i.e. EDMD \cite{Rowley_EDMD}, kernel EDMD \cite{Williams_kernel}, jetEDMD \cite{Ishikawa_JetEDMD}). For analytic EDMD and EDMD, we employ translated monomials $(\mathbf{x}-\mathbf{x}^*)^\alpha$, with $|\alpha|\leq|\alpha|_{max}$, as basis functions. For jetEDMD, we compute the EDMD matrix with translated monomials up to total degree $10$ (for the Van der Pol and Duffing systems) or up to total degree $4$ (for the $10$-dimensional networked dynamics), and then truncate the EDMD matrix up to total degree $|\alpha|_{max}$. We use the Szeg\"o kernel \mbox{$k(\mathbf{x},\mathbf{y})=\prod_{i=1}^n 1/(1-(x_i-x_i^*) (y_i-x^*_i))$} for both analytic EDMD and kernel EDMD, with no regularization. In the case of the bistable Duffing dynamics, monomials are translated with respect to the equilibrium toward which the data trajectory converges. The test points used to compute the estimated eigenfunction accuracy \cref{eq:metric3} are taken in the basin of attraction of that equilibrium.

\begin{table}[t]
    \centering
    \begin{tabular}{cccccc}
       \textbf{Method} & $\mathrm{\textbf{ESA}}_1$ & $\mathrm{\textbf{ESA}}_2$ & $\mathrm{\textbf{ESA}}_3$ & $\mathrm{\textbf{SPM}}$ & $\mathrm{\textbf{EFA}}$\\
       \hline
       \hline
       & & \textbf{Van der Pol} & ($M=75$)\\
               \hline
        analytic EDMD & $\mathbf{1.13 \times 10^{-5}}$ & $\mathbf{2.43 \times 10^{-4}}$ & $\mathbf{3.35 \times 10^{-3}}$ & $\mathbf{9.83 \times 10^{-2}}$ & $\mathbf{7.65 \times 10^{-3}}$ \\
        EDMD & $2.31 \times 10^{-2}$ & $0.23$ & $0.39$ & $4.72 \times 10^{-1}$ & $1.69 \times 10^{-2}$ \\
        kernel EDMD & $1.07 \times 10^{-4}$ & $6.80 \times 10^{-3}$ & $5.73 \times 10^{-2}$ & $2.30$ & $5.17$ \\
        jetEDMD & $4.03 \times 10^{-4}$ & $1.66 \times 10^{-2}$ & $7.51 \times 10^{-2}$ & $0.25$ & $1.26 \times 10^{-2}$\\
        \hline
        & & \textbf{Van der Pol} & ($M=250$)\\
       \hline
        analytic EDMD & $\mathbf{1.61 \times 10^{-10}}$ & $\mathbf{2.91 \times 10^{-8}}$ & $\mathbf{9.22 \times 10^{-7}}$ & $\mathbf{1.42 \times 10^{-3}}$ & $\mathbf{6.59 \times 10^{-3}}$ \\
        EDMD & $4.04 \times 10^{-2}$ & $0.21$ & $4.07 \times 10^{-1}$ & $0.47$ & $2.33 \times 10^{-2}$ \\
        kernel EDMD & $0.14$ & $2.24 \times 10^{-1}$ & $3.07 \times 10^{-1}$ & $4.44$ & $495$ \\
        jetEDMD & $1.53 \times 10^{-3}$ & $2.90 \times 10^{-2}$ & $1.54 \times 10^{-1}$ & $0.24$ & $1.73 \times 10^{-2}$ \\
        \hline
        & & \textbf{Duffing} & ($M=100$)\\
               \hline
        analytic EDMD & $\mathbf{8.14 \times 10^{-5}}$ & $\mathbf{9.61 \times 10^{-3}}$ & $\mathbf{6.83 \times 10^{-2}}$ & $\mathbf{1.48 \times 10^{-2}}$ & $\mathbf{9.81 \times 10^{-2}}$ \\
        EDMD & $6.15 \times 10^{-2}$ & $0.27$ & $0.79$ & $0.25$ & $0.10$ \\
        kernel EDMD & $8.26 \times 10^{-2}$ & $0.14$ & $0.24$ & $8.39$ & $0.23$ \\
        jetEDMD & $2.60 \times 10^{-2}$ & $0.31$ & $0.46$ & $0.68$ & $0.14$ \\
              \hline
        & & \textbf{Duffing} & ($M=250$)\\
               \hline
        analytic EDMD & $1.42 \times 10^{-7}$ & $4.27 \times 10^{-6}$ & $4.98 \times 10^{-4}$ & $1.32 \times 10^{-4}$ & $\mathbf{9.37 \times 10^{-2}}$ \\
        EDMD & $7.59 \times 10^{-2}$ & $0.28$ & $0.80$ & $0.25$ & $9.95 \times 10^{-2}$ \\
        kernel EDMD & $0.17$ & $0.25$ & $0.34$ & $10.14$ & $0.42$ \\
        jetEDMD & $\mathbf{2.07 \times 10^{-9}}$ & $\mathbf{1.90 \times 10^{-7}}$ & $\mathbf{5.02 \times 10^{-6}}$ & $\mathbf{1.50 \times 10^{-6}}$ & $\mathbf{9.37 \times 10^{-2}}$ \\
                \hline
        & & \textbf{$\mathbf{10}$-dim network} \\
               \hline
        analytic EDMD & $\mathbf{1.95 \times 10^{-3}}$ & $0.14$ & - & $\mathbf{1.15 \times 10^{-2}}$ & $0.97$ \\
        EDMD & $0.20$ & $0.26$ & - & $5.82 \times 10^{-2}$ & $\mathbf{0.43}$ \\
        kernel EDMD & $0.15$ & $0.20$ & - & $2.86 \times 10^{-2}$ & $3.61$ \\
        jetEDMD & $1.58 \times 10^{-2}$ & $\mathbf{9.28 \times 10^{-2}}$ & - & $1.68 \times 10^{-2}$ & $0.96$
    \end{tabular}
    \caption{Comparison results (averaged over $50$ simulations).}
    \label{tab:res_compar}
\end{table}

The results averaged over $50$ simulations are given in \cref{tab:res_compar}. Although this comparative study is not exhaustive, it is clear that analytic EDMD outperforms the other methods in nearly all situations. In particular, it is very accurate, even for higher orders, as shown by very small values of the ESA measure. Spectral pollution captured through the SPM measure is also very low for analytic EDMD, and the EFA measure shows that the dominant principal eigenfunction is accurately computed. Note that the SPM measure could be further improved, in particular for other methods such as kernel EDMD, by leveraging recent residual-based techniques \cite{Colbrook2024,Drmac2018}. Analytic EDMD is slightly less performant when the data set is drawn from a single trajectory (Duffing dynamics), since convergence results for the inner product approximation require a sampling over a positive measure set (see \cref{prop:bound_entries} in \cref{sec:errors} below). In particular, jetEDMD outperforms analytic EDMD in spectrum estimation when the time range of the trajectory is large (equal to $25$ in the case $M=250$), a situation where data points are available in a close neighborhood of the equilibrium. Finally, for the $10$-dimensional networked dynamics, analytic EDMD outperforms the other methods for eigenvalues estimation (except that jetEDMD outperforms for $\textrm{ESA}_2$), while EDMD provides the most accurate computation of the principal eigenfunction.

\subsection{Effect of measurement noise}

We now investigate the effect of measurement noise on the performance of analytic EDMD. To do so, $250$ data pairs are generated by the Van der Pol dynamics (see parameters in \cref{tab:simu_param}) and corrupted with additive Gaussian noise, with zero mean and pointwise standard deviation $\sigma_{noise}\in\{0.001,0.01,0.1\}$. The methods are used with the same parameters as in \cref{sec:comparison}, except that analytic EDMD and kernel EDMD are regularized (with the regularization parameter $\epsilon= \sigma_{noise}$). The results show that analytic EDMD accurately estimates both eigenvalues and eigenfunction in presence of weak measurement noise (see \cref{tab:res_noise}). As the noise level increases, one observes a performance drop, but analytic EDMD still achieves comparable performance with EDMD in the case of strong noise.

\begin{table}[t]
    \centering
    \begin{tabular}{cccccc}
       \textbf{Method} & $\mathrm{\textbf{ESA}}_1$ & $\mathrm{\textbf{ESA}}_2$ & $\mathrm{\textbf{ESA}}_3$ & $\mathrm{\textbf{SPM}}$ & $\mathrm{\textbf{EFA}}$\\
       \hline
       \hline
       & & $\sigma_{noise}=0.001$ \\
       \hline
        analytic EDMD & $\mathbf{4.13 \times 10^{-3}}$ & $\mathbf{1.60 \times 10^{-2}}$ & $\mathbf{6.97 \times 10^{-2}}$ & $\mathbf{0.21}$ & $\mathbf{1.52 \times 10^{-2}}$ \\
        EDMD & $4.61 \times 10^{-2}$ & $0.21$ & $0.43$ & $0.47$ & $2.61 \times 10^{-2}$ \\
        kernel EDMD & $5.82 \times 10^{-2}$ & $0.22$ & $0.34$ & $19.54$ & $342$ \\
        jetEDMD & $1.10 \times 10^{-2}$ & $8.46 \times 10^{-2}$ & $0.19$ & $0.34$ & $0.13$\\
       \hline
       & & $\sigma_{noise}=0.01$ \\
       \hline
        analytic EDMD & $\mathbf{2.11 \times 10^{-2}}$ & $\mathbf{8.92 \times 10^{-2}}$ & $\mathbf{0.24}$ & $\mathbf{0.28}$ & $\mathbf{4.10 \times 10^{-2}}$ \\
        EDMD & $6.41 \times 10^{-2}$ & $0.30$ & $0.48$ & $0.49$ & $4.24 \times 10^{-2}$ \\
        kernel EDMD & $0.20$ & $0.48$ & $0.48$ & $22.38$ & $257$ \\
        jetEDMD & $8.51 \times 10^{-2}$ & $0.32$ & $0.41$ & $1.02$ & $0.56$ \\
       \hline
       & & $\sigma_{noise}=0.1$ \\
       \hline
        analytic EDMD & $\mathbf{0.14}$ & $\mathbf{0.37}$ & $0.72$ & $\mathbf{0.59}$ & $0.21$ \\
        EDMD & $0.16$ & $\mathbf{0.37}$ & $\mathbf{0.54}$ & $0.75$ & $\mathbf{0.17}$ \\
        kernel EDMD & $0.24$ & $0.53$ & $0.70$ & $24.57$ & $47.90$ \\
        jetEDMD & $0.21$ & $0.49$ & $0.80$ & $2.44$ & $0.98$ 
    \end{tabular}
    \caption{Effect of measurement noise in the case of the Van der Pol dynamics. The results are averaged over $50$ simulations.}
    \label{tab:res_noise}
\end{table}

\subsection{Choice of kernel}

The performance of analytic EDMD depends on the kernel. We will investigate this effect by considering several kernels:
\begin{itemize}
    \item the Szeg\"o kernel $k(\mathbf{x},\mathbf{y})=\prod_{i=1}^n 1/(1-\gamma^2 \, x_i \, y_i)$ associated with the Hardy space over a polydisk of radius $1/\gamma$;
    \item the Szeg\"o kernel $k(\mathbf{x},\mathbf{y})=1/(1-\gamma^2 \, \mathbf{x}^T \mathbf{y})$ associated with the Hardy space over a ball of radius $1/\gamma$;
    \item the exponential kernel $k(\mathbf{x},\mathbf{y})=\exp(\gamma^2 \, \mathbf{x}^T \mathbf{y})$ associated with the Fock space of analytic entire functions;
    \item the polynomial kernel $k(\mathbf{x},\mathbf{y})=(1+\gamma^2 \, \mathbf{x}^T \mathbf{y})^d$, with $d=20$, associated with a finite-dimensional space spanned by monomials.\\
\end{itemize}

\begin{figure}[h!]
\centering
\subfigure[Estimated spectrum accuracy]{\includegraphics[width=0.49\linewidth]{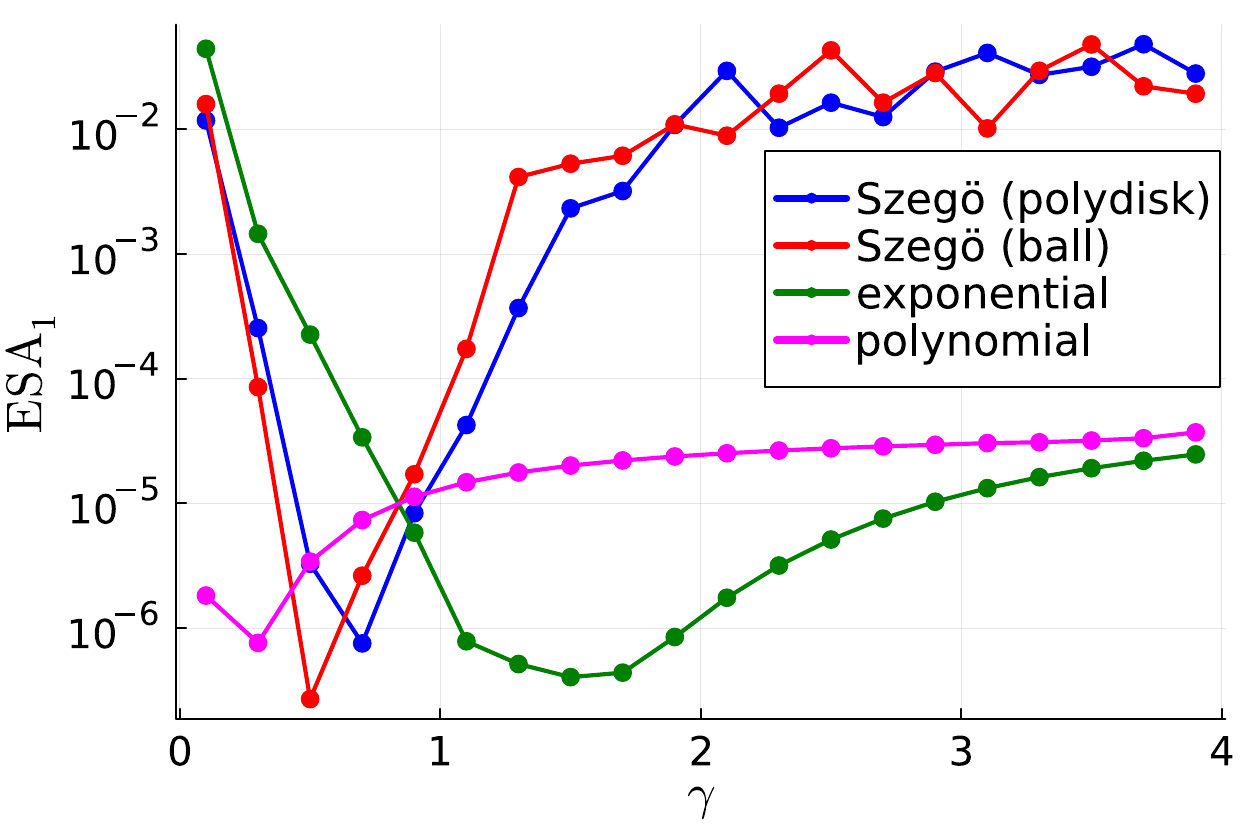}}
\subfigure[Spectral pollution measure]{\includegraphics[width=0.49\linewidth]{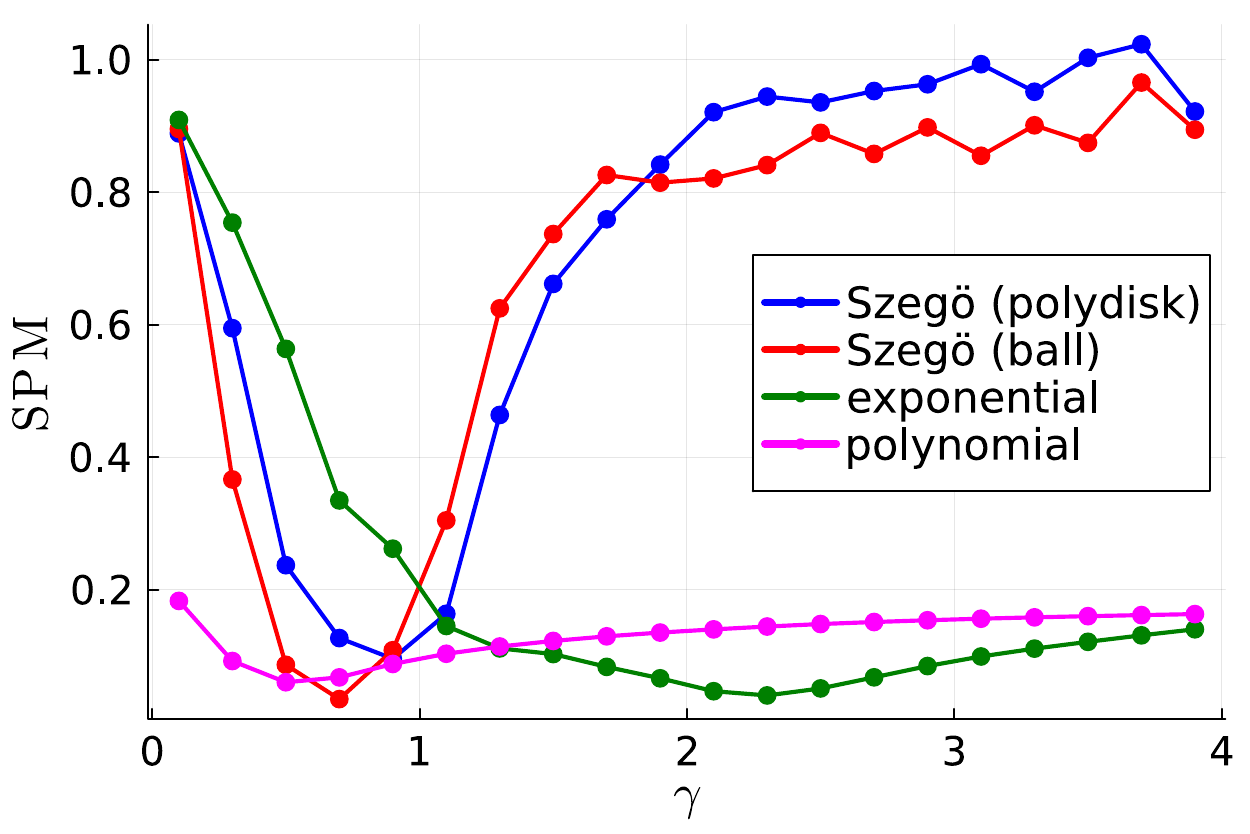}}
\caption{Effect of the kernel on the estimated spectrum accuracy $\textrm{ESA}_1$ (first order) and the spectral pollution measure $\textrm{SPM}$. The results are averaged over $50$ simulations.}
\label{fig_kernels}
\end{figure}

Note that the hyperparameter $\gamma>0$ acts as a scaling factor that should be adapted to the sampling region. All the kernels are of Taylor-type (except the first one) and yield an inner product that makes the monomials orthogonal. We use \cref{eq:approx_non_ortho} with a basis of unweighted monomials to compute the Koopman matrix and focus on spectrum estimation. The results shown in \cref{fig_kernels} are obtained for the Van der Pol dynamics with the parameters given in \cref{tab:simu_param} in the case $M=75$. We observe that all kernels provide comparatively good results for well-chosen values of the scaling factor $\gamma$, which either allows to consider a larger domain ($\gamma<1$ for the Szeg\"o and polynomial kernels) or a smaller one ($\gamma>1$ for the exponential kernel). The choice of scaling factor is particularly critical for the Szeg\"o kernel (see e.g. \cref{fig_kernels}(b)), due to singularities lying on the polydisk or ball of radius $1/\gamma$. Optimizing the kernels through a proper tuning of the hyperparameters is beyond the scope of this paper but could potentially be achieved by leveraging the expected lattice structure of the spectrum. 

\section{Error bounds and convergence}
\label{sec:errors}

In this section, we compute upper bounds on the estimation error of the Koopman operator spectrum, and establish the convergence of both eigenvalues and eigenfunctions approximations. Note that spectral bounds have also been proposed in the context of kernel-based learning of stochastic Koopman operators in \cite{Kostic_spectral_bounds}.

The following result provides an upper estimate of the error between the entries $\mathbf{K}_{ij}$ of the Koopman matrix obtained with the exact Taylor projection and the entries $\hat{\mathbf{K}}_{ij}$  of the matrix \cref{eq:analytic_EDMD} computed from the data.
\begin{proposition}
\label{prop:bound_entries}
    Let $\hat{\mathbf{K}}_{ij}= \mathbf{e}_i^T  \mathbf{G}^{-1} \mathbf{e}_j'$ and $\mathbf{K}_{ij}=\langle e_i,K e_j \rangle_\mathcal{H}$. Then,
    \begin{equation}
    \label{eq:bound_entries}
|\mathbf{K}_{ij} - \mathbf{\hat{K}}_{ij}| \leq \sqrt{1-\mathbf{e}^T_i \mathbf{G}^{-1} \mathbf{e}_i} \, \|K e_j\|_\Hc.
\end{equation}
Moreover\footnote{Part of this result was initially mentioned by Isao Ishikawa in a private communication.}, if the data points $\{\mathbf{x}_k\}_{k=1}^M$ are randomly drawn from a continuous distribution with support $\bar{X} \subseteq X$ of positive Lebesgue measure, then
\begin{equation}
\label{eq:conv_entries}
    \lim_{M \rightarrow \infty} |\mathbf{K}_{ij} - \mathbf{\hat{K}}_{ij}| = 0 \textrm{ almost surely.}
\end{equation}
\end{proposition}
\begin{proof}
It follows from \cref{rem:inner_product} that
\begin{equation}
\label{eq:diff_Kij}
\hat{\mathbf{K}}_{ij} - \mathbf{K}_{ij} = \langle \Pi e_i,K e_j \rangle_\mathcal{H} - \langle e_i,K e_j \rangle_\mathcal{H} = \langle \Pi e_i - e_i, K e_j \rangle_\mathcal{H}.
\end{equation}
Then, using the Cauchy-Schwarz inequality and basic properties of the orthogonal projection, we have
\begin{equation*}
|\mathbf{K}_{ij} - \mathbf{\hat{K}}_{ij}| \leq \|\Pi e_i - e_i\|_\Hc \, \|K e_j\|_\Hc = \sqrt{\|e_i\|^2_\Hc-\|\Pi e_i\|_\Hc^2} \, \|K e_j\|_\Hc  .
\end{equation*}
The inequality \cref{eq:bound_entries} is obtained by observing that $\|e_i\|^2_\Hc=1$ and
$$\|\Pi e_i\|^2_\Hc = \langle \Pi e_i,\Pi e_i\rangle_\Hc = \langle \Pi e_i, e_i\rangle_\Hc=\mathbf{e}^T_i \mathbf{G}^{-1} \mathbf{e}_i,$$
where we used \cref{eq:prop_inner_prod}.

Finally, suppose that an analytic function $f \in \Hc$ satisfies $0=\langle f, k_{x_j} \rangle=f(x_j)$ for all $j$. Since the set $\{x_j\}_{j=1}^\infty$ is (almost surely) dense in $\bar{X}$, this implies that $f(x)=0$ for all $x\in \bar{X}$ by continuity of $f$ and therefore $f=0$ according to the identity theorem \cite{identity_theorem}. Hence, the linear span of $\{k_{x_j}\}_{j=1}^\infty$ is dense in $\Hc$, so that the orthogonal projection $\Pi$ satisfies
\begin{equation*}
    \lim_{M \rightarrow \infty} \|\Pi e_i-e_i\|_\Hc = 0
\end{equation*}
and the result follows from \cref{eq:diff_Kij}.
\end{proof}
It is noticeable that convergence holds if the data points are sampled according to a general distribution and from an arbitrary set (of positive measure) of the state space. This confirms the numerical result shown in \cref{subsec:specific_cases}.

Next, we provide a bound on the estimation error of the Koopman eigenvalues.
\begin{proposition}
\label{prop:bounds_eigen}
Let $\mu_j$ be the eigenvalues of the Jacobian matrix of the map $\varphi$ evaluated at the equilibrium.
    Let $r\in\{1,\dots,|\alpha|_{max}\}$ and define $l_{min}$ and $l_{max}$ as, respectively, the minimum and maximum row/column indices of the matrix block $\overline{\mathbf{K}}_{rr}$, i.e. $l_{min}=\min \{j:|\overline{\alpha}(j)|=r\}$ and $l_{max}=\max \{j:|\overline{\alpha}(j)|=r\}$. If $\overline{\mathbf{K}}_{rr}$ is diagonalizable, then for any Koopman eigenvalue $\mu_1^{\alpha_1} \cdots \mu_n^{\alpha_n}$, with $|\alpha|=r$, there exists an eigenvalue $\hat{\mu}\in \sigma(\overline{\mathbf{K}}_{rr}) \subset \sigma(\hat{\mathbf{K}})$, such that
    \begin{equation*}
        \left|\hat{\mu}-\mu_1^{\alpha_1} \cdots \mu_n^{\alpha_n}\right| \leq \kappa_1(\mathbf{V}) \, \max_{l_{min}\leq j \leq l_{max}} \|K e_j\|_\Hc \, \sum_{i=l_{min}}^{l_{max}} \sqrt{1-\mathbf{e}_i \mathbf{G}^{-1} \mathbf{e}_i},
    \end{equation*}
    \begin{equation*}
        \left|\hat{\mu}-\mu_1^{\alpha_1} \cdots \mu_n^{\alpha_n}\right| \leq \kappa_2(\mathbf{V}) \, \sqrt{\sum_{i=l_{min}}^{l_{max}} \sum_{j=l_{min}}^{l_{max}}  \left( 1-\mathbf{e}_i \mathbf{G}^{-1} \mathbf{e}_i \right) \, \|K e_j\|^2_\Hc},
    \end{equation*}  
    and
    \begin{equation*}
        \left|\hat{\mu}-\mu_1^{\alpha_1} \cdots \mu_n^{\alpha_n}\right| \leq \kappa_\infty(\mathbf{V}) \, \max_{l_{min}\leq i \leq l_{max}} \sqrt{1-\mathbf{e}_i \mathbf{G}^{-1} \mathbf{e}_i} \, \sum_{j=l_{min}}^{l_{max}} \|K e_j\|_\Hc,
    \end{equation*}
where $\mathbf{V}$ is the matrix of eigenvectors of $\overline{\mathbf{K}}_{rr}$ and $\kappa_p(\mathbf{V}) = \|\mathbf{V}\|_p \|\mathbf{V}^{-1}\|_p $, with $p\in\{1,2,\infty\}$, is the condition number of $\mathbf{V}$. 
\end{proposition}
\begin{proof}
    Consider the matrix block $\mathbf{\tilde{K}}_{rr}$ of $\mathbf{K}$ corresponding to the same rows and columns as $\overline{\mathbf{K}}_{rr}$ in $\hat{\mathbf{K}}$. It is known that the eigenvalues of $\mathbf{\tilde{K}}_{rr}$ are Koopman eigenvalues $\mu_1^{\alpha_1} \cdots \mu_n^{\alpha_n}$, with $|\alpha|=r$. Moreover, it follows from the Bauer-Fike theorem \cite[Theorem III.a]{Bauer_Fike} that, for all  $\mu_1^{\alpha_1} \cdots \mu_n^{\alpha_n} \in \sigma(\mathbf{\tilde{K}}_{rr})$, there exists $\hat{\mu} \in \sigma(\overline{\mathbf{K}}_{rr})$ such that
    \begin{equation*}
        \left|\hat{\mu}-\mu_1^{\alpha_1} \cdots \mu_n^{\alpha_n}\right| \leq \kappa_p(\mathbf{V}) \, \|\overline{\mathbf{K}}_{rr}-\mathbf{\tilde{K}}_{rr}\|_p.
    \end{equation*}
    Then, the result follows from \cref{eq:bound_entries} in \cref{prop:bound_entries}, the definition of the matrix $1$ and $\infty$-norms, and the fact that the matrix $2$-norm is bounded by the Frobenius norm. This concludes the proof.
\end{proof}
This proposition provides a bound for the estimated spectrum accuracy \cref{eq:metric1}. The argument based on the Bauer-Fike theorem could also be ``reversed'' to provide a bound for the spectral pollution measure \cref{eq:metric2}, where $\mathbf{V}$ is replaced by the matrix of eigenvectors of $\mathbf{\tilde{K}}_{rr}$. Moreover, we note that the error bounds only depend on the set of data points and on the map $\varphi$, through the quantities $\sqrt{1-\mathbf{e}_i \mathbf{G}^{-1} \mathbf{e}_i}$ and $\|K e_j\|_\Hc$, respectively. The former can be computed from the data, without knowledge on the vector field, and therefore can be used to choose the best data set minimizing the error. The quantity $\|K e_j\|_\Hc$ can possibly be obtained from some prior knowledge on the map $\varphi$. For instance, in the Hardy space on the polydisk, we have $\|K e_j\|_\Hc = \|\varphi^{\overline{\alpha}(j)}\|_\Hc \leq (\max_{z \in \partial \mathbb{D}^n} \|\phi(z))\|_\infty )^{|\overline{\alpha}(j)|} \triangleq \varphi_{max}^{|\overline{\alpha}(j)|}$ where $\partial \mathbb{D}^n$ is the boundary of the unit polydisk.

\begin{remark}[Continuous time]
    In the continuous-time case, the Koopman eigenvalues of the form $\alpha_1 \lambda_1 + \dots + \alpha_n \lambda_n$ are estimated by $\hat{\lambda} = \log(\hat{\mu})/\Delta t$, so that the error on $\hat{\mu}$ propagates to $\hat{\lambda}$ through the complex logarithm. Assume that we have a bound $\delta$ such that $\left|\hat{\mu}-\mu_1^{\alpha_1} \cdots \mu_n^{\alpha_n}\right|<\delta$ (see above) and denote $\nabla \log(x+iy) = \left(\frac{\partial \log}{\partial x}(x+iy), \frac{\partial \log}{\partial y}(x+iy)\right)$. Straightforward computations yield $\|\nabla \log(z)\|=\sqrt{2}/|z|$ and we have
    \begin{equation*}
    \begin{split}
        |\hat{\lambda}-(\alpha_1 \lambda_1 + \dots + \alpha_n \lambda_n)| = \frac{1}{\Delta t} \left|\log(\hat{\mu})-\log\left(\mu_1^{\alpha_1} \cdots \mu_n^{\alpha_n}\right)\right| & \leq \frac{\delta}{\Delta t} \max_{|z-\hat{\mu}|<\delta} \|\nabla \log(z)\| \\
        & \leq \frac{\sqrt{2}}{\Delta t} \frac{\delta}{|\hat{\mu}|-\delta}
        \end{split}
    \end{equation*}
     provided that $|\hat{\mu}|>\delta$. \hfill $\diamond$
\end{remark}

Finally, the following proposition ensures convergence for both Koopman eigenvalues and eigenfunctions estimations.
\begin{proposition}
\label{prop:convergence}
If the data points $\{\mathbf{x}_k\}_{k=1}^M$ are randomly drawn from a continuous distribution with support $\bar{X} \subseteq X$ of positive Lebesgue measure, then for any Koopman eigenvalue $\mu_1^{\alpha_1} \cdots \mu_n^{\alpha_n}$, with $|\alpha|=r$, there exists a sequence of eigenvalues $\hat{\mu}^{(M)}$ of $\overline{\mathbf{K}}_{rr}$ such that
\begin{equation}
\label{eq:conv_eigenval}
    \lim_{M \rightarrow \infty} \left|\hat{\mu}^{(M)}-\mu_1^{\alpha_1} \cdots \mu_n^{\alpha_n}\right| = 0 \quad \textrm{almost surely}.
\end{equation}
Moreover, suppose that the Jacobian eigenvalues $\mu_j$ are distinct. Then the vectors $\overline{\mathbf{v}}_r^{(j)}$ obtained in \cref{eq:eigvec} converge to the $r$-th order Taylor coefficients of the eigenfunction $\phi_{\mu_j}$ (in the basis $\{e_i\}_{i=1}^N$) and the associated approximated eigenfunction $\hat{\phi}_{\mu_j}^{(M,N)}$ computed with \cref{eq:taylor_Koop_eigfct} satisfies
\begin{equation*}
    \lim_{N \rightarrow \infty} \lim_{M \rightarrow \infty} \|\hat{\phi}_{\mu_j}^{(M,N)}-\phi_{\mu_j}\|_{\Hc} = 0 \quad \textrm{almost surely}.
\end{equation*}
\end{proposition}
\begin{proof}
We have that \cref{eq:conv_entries} in \cref{prop:bound_entries} implies that $\overline{\mathbf{K}}_{rr} \rightarrow \mathbf{\tilde{K}}_{rr}$ as $M \rightarrow \infty$ and recall that $\mu_1^{\alpha_1} \cdots \mu_n^{\alpha_n} \in \sigma(\mathbf{\tilde{K}}_{rr})$. Then, \cref{eq:conv_eigenval} follows since the eigenvalues depend continuously on the matrix entries.

Since an eigenvector associated with a distinct eigenvalue depends continuously on the matrix entries, we have that $\overline{\mathbf{v}}_1^{(j)} \rightarrow \tilde{\mathbf{v}}_1^{(j)}$, where $\overline{\mathbf{v}}_1^{(j)}$ and $\tilde{\mathbf{v}}_1^{(j)}$ are the (unit) eigenvectors of the matrices $\overline{\mathbf{K}}_{11}$ and $\mathbf{\tilde{K}}_{11}$, respectively. According to \cref{eq:eigvec}, we have
\begin{equation*}
    \overline{\mathbf{v}}_r^{(j)} = (\overline{\mathbf{K}}_{rr} - \hat{\mu}_j^{(M)} I)^{-1} \sum_{s=1}^{r-1} \overline{\mathbf{K}}_{rs} \overline{\mathbf{v}}_s^{(j)} \qquad
    \tilde{\mathbf{v}}_r^{(j)} = (\tilde{\mathbf{K}}_{rr} - \mu_j I)^{-1} \sum_{s=1}^{r-1} \tilde{\mathbf{K}}_{rs} \tilde{\mathbf{v}}_s^{(j)}
\end{equation*}
for $r>1$, where $\tilde{\mathbf{v}}_r^{(j)}$ denotes the vector of exact Taylor coefficients. Given that $\overline{\mathbf{K}}_{rs} \rightarrow \mathbf{\tilde{K}}_{rs}$ for all $r,s$ and $\hat{\mu}_j^{(M)} \rightarrow \mu_j$ as $M \rightarrow \infty$ (\cref{prop:bound_entries}), a recursive argument directly shows that $\overline{\mathbf{v}}_r^{(j)} \rightarrow \tilde{\mathbf{v}}_r^{(j)}$ for all $r$ as $M \rightarrow \infty$. Finally, we have
\begin{equation*}
    \|\hat{\phi}_{\mu_j}^{(M,N)}-\phi_{\mu_j}\|_{\Hc} \leq \left\|\sum_{r=1}^{|\alpha|_{max}} (\overline{\mathbf{v}}_r^{(j)}-\tilde{\mathbf{v}}_r^{(j)})^T \, \mathbf{\overline{e}}_r \right\|_{\Hc} + \left\|\sum_{r=1}^{|\alpha|_{max}} (\tilde{\mathbf{v}}_r^{(j)})^T \, \mathbf{\overline{e}}_r - \phi_{\mu_j} \right\|_{\Hc}.
\end{equation*}
The first term of the right-hand side converges to $0$ as $M \rightarrow \infty$, while the second term does not depend on $M$ and converges to $0$ as $N \rightarrow \infty$ since the monomials $\{e_i\}_{i=1}^\infty$ form a basis in $\Hc$. This concludes the proof.
\end{proof}
The above results provide the convergence of a sequence of estimated eigenvalues and strong convergence in the RKHS norm of the estimated eigenfunctions, a property that contrasts with the weaker convergence properties of EDMD (i.e. convergence of a subsequence of estimated eigenvalues and weak convergence of estimated eigenfunctions, see \cite{Korda_convergence}). It is also remarkable that the convergence of the eigenvalues does not depend on the maximal degree of the monomials basis. In particular, the Jacobian eigenvalues can be estimated with arbitrary accuracy by using only monomials of degree $1$.

\begin{example}
We illustrate the results of \cref{prop:bounds_eigen} and \cref{prop:convergence} with the map
\begin{equation}
\label{eq:disc_map}
    \varphi(x_1,x_2) = (0.2\, x_1-0.5\, x_1\,x_2, 0.3\, x_2 + 0.6 \, x_1\, x_2),
\end{equation}
which possesses a stable equilibrium at the origin. The data points are distributed over $[0,1]^2$. Over $50$ simulations, we compute the maximal error on the eigenvalues of order $1$ and $2$, which corresponds to the estimated spectrum accuracy $\textrm{ESA}_1$ and $\textrm{ESA}_2$, respectively. Error bounds are obtained according to \cref{prop:bounds_eigen} with $\|K e_j\|_\Hc \leq \varphi_{max}^{|\overline{\alpha}(j)|}$ and with the prior knowledge $\varphi_{max}=0.8$. Since the error bound depends on the data set, it varies from one simulation to another but is verified to always be greater than the corresponding true error. The results shown in \cref{fig_error} are consistent with the convergence properties for large $M$, and demonstrate the high accuracy of the method, even for a low number of data points (see e.g. $M=30$). Moreover, for first order eigenvalues (i.e. Jacobian eigenvalues), the error bounds guarantee an error less than $10^{-2}$ in average with $M=50$, and less than $10^{-3}$ in average with $M=100$ data pairs, while the average true error is of the order of $10^{-4}$ and $10^{-7}$, respectively.

\begin{figure}[h!]
\centering
\subfigure[Order $|\alpha|=1$]{\includegraphics[width=0.49\linewidth]{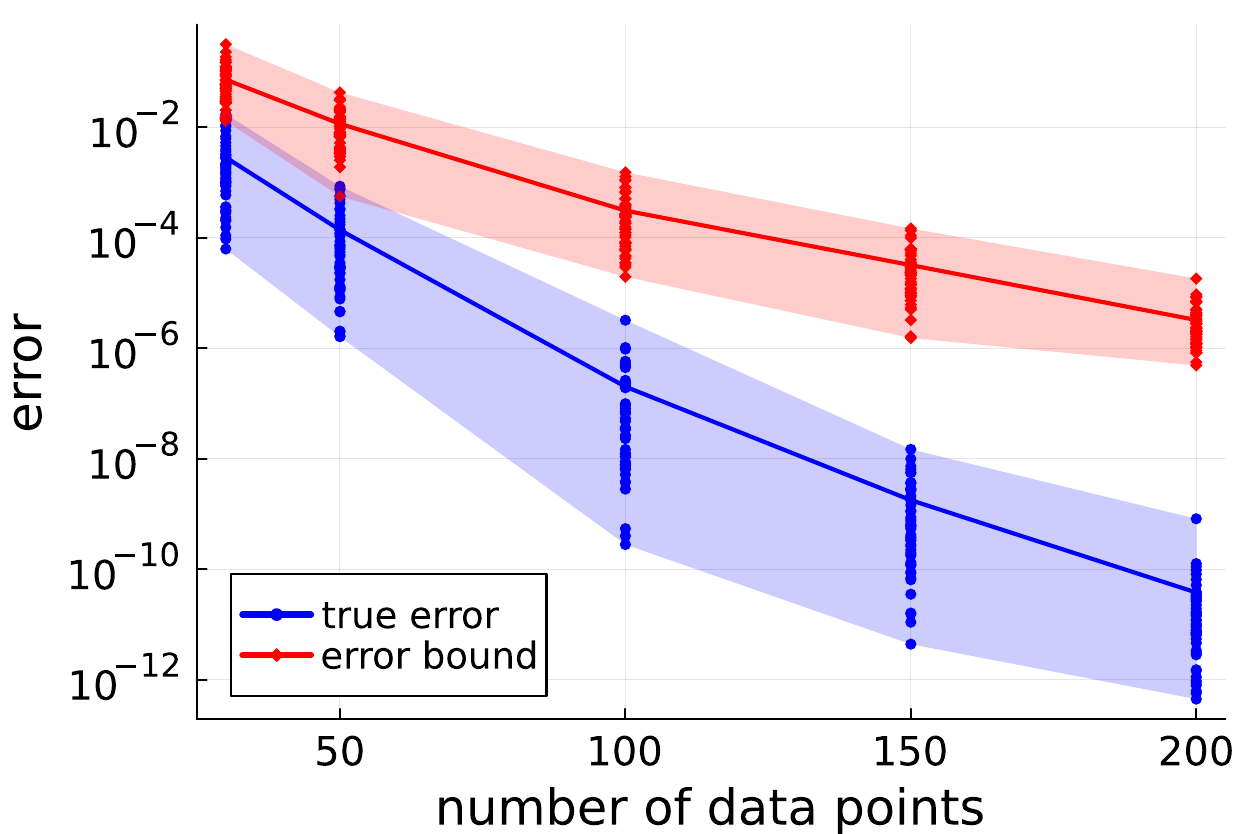}}
\subfigure[Order $|\alpha|=2$]{\includegraphics[width=0.49\linewidth]{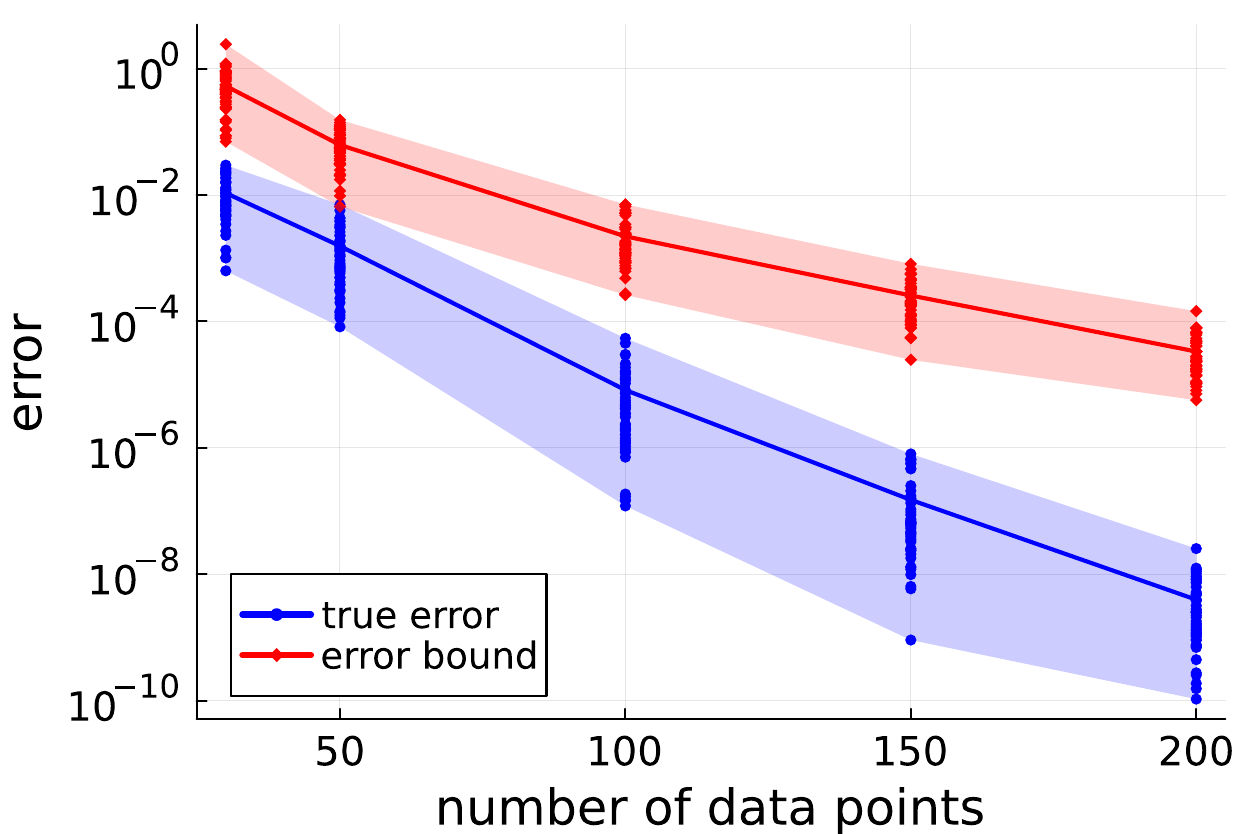}}
\caption{True error (estimated spectrum accuracy $\textrm{ESA}_\alpha$ with (a) $\alpha=1$ and (b) $\alpha=2$) and error bounds on the computation of Koopman eigenvalues for the map \cref{eq:disc_map}. Each dot corresponds to one simulation. The solid lines represent the errors averaged over $50$ simulations.}
\label{fig_error}
\end{figure}

\end{example}

\section{Conclusion}
\label{sec:conclu}

Analytic EDMD is a simple EDMD-type method that reveals the spectral properties of the Koopman operator in spaces of analytic functions. The method relies on orthogonal Taylor projections in RKHS onto a subspace spanned by an orthonormal basis of monomial functions. It is shown to yield excellent performance to capture the lattice-structured spectrum of the Koopman operator and the corresponding analytic eigenfunctions in the case of dynamics with hyperbolic equilibria. In particular, the method does not suffer from spectral pollution and can reach arbitrary accuracy on the spectrum with a fixed finite dimension of the approximation, even if no data points are taken from the neighborhood of the equilibrium. In addition, analytic EDMD is complemented with theoretical convergence results and upper bounds on the spectrum estimation error.

We envision several research perspectives that could complement this work. First, analytic EDMD could be extended to limit-cycle dynamics which are also associated with a block-triangular structure of the Koopman operator, or could be adapted to high-dimensional dynamics by using the kernel trick. Moreover, its performance could be further investigated through tighter error bounds, convergence rate characterization, and noise robustness analysis. On the same line, future work could provide prediction error bounds similar to those obtained in \cite{Worthmann_kernel_errors}. The method could also be combined with optimization or machine learning techniques for a proper choice of kernel and tuning of kernel parameters. Finally, one could investigate the potential of analytic EDMD to improve existing Koopman operator-based methods in various contexts such as control theory, stability analysis, parameter estimation, isostable reduction, and spectral network identification.

\section*{Acknowledgments}

This work was initiated while A. Mauroy was a visiting researcher in I. Mezi\'c's group at UCSB, supported by a grant from the Belgian National Funds for Scientific Research (F.R.S.-FNRS). A. Mauroy thanks F.-G. Bierwart and C. Mugisho for fruitful discussions on the triangular structure of the Koopman operator, N. Atanasov for directing his attention to Mercer's theorem, and P. Bevanda and Z. Drma\v{c} for useful discussion on kernel-based methods. The authors are also thankful to I. Ishikawa for pointing out the inner product interpretation given in \cref{rem:inner_product} and used in the theoretical results of \cref{sec:errors}.

I. Mezi\'c's work was partially supported by AFOSR award number FA9550-22-1-0531 and also based upon work supported by the Defense Advanced Research Projects Agency (DARPA) under Agreement No. HR00112590152. Approved for public release; distribution is unlimited.

This research used resources of the ``Plateforme Technologique de Calcul Intensif (PTCI)'' located at the University of Namur, Belgium, which is supported  by the FNRS-FRFC, the Walloon Region, and the University of Namur (Conventions No. 2.5020.11, GEQ U.G006.15, 1610468, RW/GEQ2016 et U.G011.22). The PTCI is member of the ``Consortium des Equi-pements de Calcul Intensif  (CECI)''.

\bibliographystyle{siamplain}

\end{document}